\documentclass[a4paper,12pt,reqno]{amsart}
\usepackage[T1]{fontenc}
\usepackage[utf8]{inputenc}
\usepackage[english]{babel}
\usepackage{amsmath,amsthm,amssymb}
\usepackage{fullpage}
\usepackage{xcolor}
\usepackage{siunitx}
\usepackage{tikz}
\usepackage{multirow}
\usepackage{pgfplotstable}
\usetikzlibrary{
	matrix,
}
\pgfplotsset{compat = newest, tick label style = {font = \tiny}}
\usepackage[width=0.95\textwidth]{caption}
\usepackage[width=0.95\textwidth,labelfont=rm]{subcaption}
\usepackage{doi,hyperref}
\makeatletter
\def\@seccntformat#1{%
  \protect\textup{\protect\@secnumfont
    \ifnum\pdfstrcmp{subsection}{#1}=0 \bfseries\fi
    \csname the#1\endcsname
    \protect\@secnumpunct
  }%
}
\makeatother
\newtheorem{theorem}{Theorem}[section]
\newtheorem{proposition}[theorem]{Proposition}
\newtheorem{lemma}[theorem]{Lemma}
\newtheorem{algorithm}[theorem]{Algorithm}
\newtheorem{definition}[theorem]{Definition}
\newtheorem{remark}[theorem]{Remark}
\newtheorem{corollary}[theorem]{Corollary}
\renewcommand{\aa}{\boldsymbol{a}}
\newcommand{\bb}{\boldsymbol{b}}
\newcommand{\cc}{\boldsymbol{c}}
\newcommand{\ee}{\boldsymbol{e}}
\newcommand\hh{\boldsymbol{h}}
\newcommand\mm{\boldsymbol{m}}

\newcommand\vv{\boldsymbol{v}}

\newcommand\CC{\boldsymbol{C}}
\newcommand\DD{\boldsymbol{D}}
\newcommand\HH{\boldsymbol{H}}
\newcommand\LL{\boldsymbol{L}}
\newcommand\MM{\boldsymbol{M}}
\newcommand\E{\mathcal{E}}
\newcommand\NN{\mathcal{N}}

\renewcommand\S{\mathcal{S}}
\newcommand\T{\mathcal{T}}
\newcommand\KKK{\boldsymbol{\mathcal{K}}}
\newcommand\RRR{\mathcal{R}}
\newcommand\XXX{\mathcal{X}}
\newcommand\HHH{\boldsymbol{\mathcal{H}}}
\newcommand\cH{c_{\HHH}}
\newcommand\cT{c_{\T}}
\newcommand\pphi{\boldsymbol{\phi}}
\newcommand\nnu{\boldsymbol{\nu}}
\newcommand\ppsi{\boldsymbol{\psi}}
\newcommand\eps{\varepsilon}
\newcommand\vvphi{\boldsymbol{\varphi}}
\newcommand\interp{\mathcal{I}_h}
\newcommand\0{\boldsymbol{0}}
\newcommand\sphere{\mathbb{S}^2}
\newcommand\curl{\nabla\times}
\renewcommand\div{\nabla\cdot}
\newcommand\Grad{\boldsymbol{\nabla}}
\newcommand\Lapl{\boldsymbol{\Delta}}
\newcommand\R{\mathbb{R}}
\newcommand\N{\mathbb{N}}
\newcommand{\abs}[1]{\left\lvert #1 \right\rvert}
\newcommand{\dual}[3][]{\langle #2,#3 \rangle_{#1}}
\newcommand{\inner}[3][]{\langle #2,#3 \rangle_{#1}}
\newcommand{\norm}[2][]{\lVert #2 \rVert_{#1}}
\DeclareMathOperator{\diam}{diam}

\newcommand\ddt{\frac{\mathrm{d}}{\mathrm{d}t}}
\newcommand\dt{\mathrm{d}t}

\newcommand\mmt{\partial_t \mm}
\newcommand\Heff[1]{{\HH_{\mathrm{eff},#1}}}
\newcommand\Hext{\HH_{\mathrm{ext}}}
\newcommand\Hstray{\HH_{\mathrm{s}}}
\newcommand\heff[1]{{\hh_{\mathrm{eff},#1}}}
\newcommand\hext{\hh_{\mathrm{ext}}}
\newcommand\Ms{M_{\mathrm{s}}} 
\newcommand\Msl[1]{M_{\mathrm{s},#1}} 

\newcommand{\weakstarto}{\overset{\ast}{\rightharpoonup}}
\newcommand{\weakto}{\rightharpoonup}

\begin{document}
\title{Convergent finite element methods\\ for antiferromagnetic and ferrimagnetic materials}
\author{Hywel~Normington}
\address{Department of Mathematics and Statistics,
University of Strathclyde,
26 Richmond Street, Glasgow G1 1XH, United Kingdom}
\email{hywel.normington@strath.ac.uk}
\author{Michele~Ruggeri}
\address{Department of Mathematics, University of Bologna, Piazza di Porta San Donato 5, 40126 Bologna, Italy}
\email{m.ruggeri@unibo.it}
\date{\today}
\keywords{antiferromagnetism; ferrimagnetism; finite element method;
$\Gamma$-convergence; Landau--Lifshitz--Gilbert equation}
\subjclass[2010]{35K61, 65M12, 65M60, 65Z05}

\begin{abstract}
We consider the numerical approximation of a continuum model of antiferromagnetic and ferrimagnetic materials.
The state of the material is described in terms of two unit-length vector fields,
which can be interpreted as the magnetizations averaging the spins of two sublattices.
For the static setting, which requires the solution of a constrained energy minimization problem,
we introduce a discretization based on first-order finite elements and prove its $\Gamma$-convergence.
Then, we propose and analyze two iterative algorithms
for the computation of low-energy stationary points.
The algorithms are obtained from (semi-)implicit time discretizations of gradient flows of the energy.
Finally, we extend the algorithms to the dynamic setting,
which consists of a nonlinear system of two Landau--Lifshitz--Gilbert equations solved by the two fields,
and we prove unconditional stability and convergence
of the finite element approximations toward a weak solution of the problem.
Numerical experiments assess the performance of the algorithms and
demonstrate their applicability
for the simulation of physical processes involving antiferromagnetic and ferrimagnetic materials.
\end{abstract}

\maketitle

\section{Introduction}

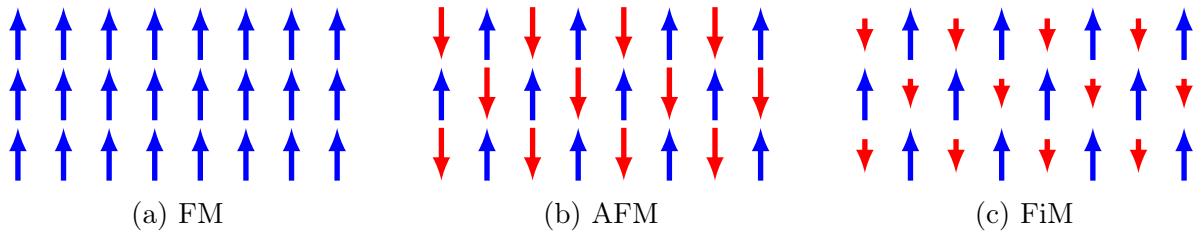
\begin{figure}[h]
\centering
\begin{subfigure}[b]{0.3\textwidth}
\centering
\begin{tikzpicture}
\foreach \x in {0,1,...,7}{\foreach \y in {0,1,2}{\draw [-latex,blue,line width=2pt](0.6*\x,0.8*\y+0.15) -- (0.6*\x,0.8*\y+0.85);}}
\end{tikzpicture}
\caption{FM}
\label{fig:fm}
\end{subfigure}
\hfill
\begin{subfigure}[b]{0.3\textwidth}
\centering
\begin{tikzpicture}
\foreach \x in {0,2,4,6}{\foreach \y in {2}{\draw [-latex,blue,line width=2pt](0.6*\x,0.8*\y+0.15) -- (0.6*\x,0.8*\y+0.85);}}
\foreach \x in {1,3,5,7}{\foreach \y in {1,3}{\draw [-latex,blue,line width=2pt](0.6*\x,0.8*\y+0.15) -- (0.6*\x,0.8*\y+0.85);}}
\foreach \x in {1,3,5,7}{\foreach \y in {2}{\draw [-latex,red,line width=2pt](0.6*\x,0.8*\y+0.85) -- (0.6*\x,0.8*\y+0.15);}}
\foreach \x in {0,2,4,6}{\foreach \y in {1,3}{\draw [-latex,red,line width=2pt](0.6*\x,0.8*\y+0.85) -- (0.6*\x,0.8*\y+0.15);}}
\end{tikzpicture}
\caption{AFM}
\label{fig:afm}
\end{subfigure}
\hfill
\begin{subfigure}[b]{0.3\textwidth}
\centering
\begin{tikzpicture}
\foreach \x in {0,2,4,6}{\foreach \y in {2}{\draw [-latex,blue,line width=2pt](0.6*\x,0.8*\y+0.15) -- (0.6*\x,0.8*\y+0.85);}}
\foreach \x in {1,3,5,7}{\foreach \y in {1,3}{\draw [-latex,blue,line width=2pt](0.6*\x,0.8*\y+0.15) -- (0.6*\x,0.8*\y+0.85);}}
\foreach \x in {1,3,5,7}{\foreach \y in {2}{\draw [-latex,red,line width=2pt](0.6*\x,0.8*\y+0.7) -- (0.6*\x,0.8*\y+0.3);}}
\foreach \x in {0,2,4,6}{\foreach \y in {1,3}{\draw [-latex,red,line width=2pt](0.6*\x,0.8*\y+0.7) -- (0.6*\x,0.8*\y+0.25);}}
\end{tikzpicture}
\caption{FiM}
\label{fig:fim}
\end{subfigure}
\caption{Classes of magnetic materials.}
\label{fig:magnetic_materials}
\end{figure}

Antiferromagnetic (AFM) and ferrimagnetic (FiM) materials,
materials in which neighboring magnetic moments
tend to align antiparallel to each other (see Figure~\ref{fig:magnetic_materials}),
have been known for many years.
However, they have recently gained renewed interest,
because several theoretical and
experimental studies have shown that AFM and FiM materials have features that could lead to strong
improvements of the functionality of spintronics devices,
compared to those based on ferromagnetic (FM) materials~\cite{bmtmot2018,kblkory2022}.

In this work,
we consider a continuum model that is the state of the art for micromagnetic simulations
of devices based on magnetic processes involving of AFM and FiM materials;
see, e.g., the works~\cite{ne2015,pkcatsf2019,spkcf2020,tslggcaf2020,sztfm2020} on AFM materials
and~\cite{mlp2016,mra2019,cstcf2021} on FiM materials.
The main elements of the model, an extension of the classical micromagnetic model of FM materials~\cite{brown1963},
are an order parameter, which consists of two unit-length vector fields
that can be interpreted as the normalized magnetizations averaging the magnetic moments of two sublattices,
and an energy functional, which consists of several contributions, each of them representing a specific physical effect.
A key feature of the model, that is necessary to describe the antiparallel alignment of the spins in AFM and FiM materials,
is a more complex expression of the exchange energy than in FM materials,
which involves not only the classical Heisenberg exchange interaction penalizing nonuniform configurations,
but also the interaction of the two fields with each other
(see the last two terms in the energy functional~\eqref{eq:energy} below). 
Similarly to the classical micromagnetic theory,
the static problem consists of minimizing the energy functional over all pairs of unit-length vector fields,
whereas the dynamics of each field out of equilibrium is governed by the Landau--Lifshitz--Gilbert (LLG) equation
(see~\eqref{eq:llg} below),
with the effective field being the Gateaux derivative of the energy with respect to the respective field.
However, due to the energy contributions involving the interaction between the two fields,
both the Euler--Lagrange equations associated with the minimization problem
and the system of LLG equations are nonlinearly coupled systems of nonlinear partial differential equations.

Building on previous work on the approximation of (the heat flow of) harmonic maps~\cite{bartels2016}
and of the classical model of FM materials~\cite{alouges2008a,ahpprs2014},
we propose fully discrete numerical schemes for the approximation of
both the static and the dynamic problems.

For the static problem,
we propose a discretization based on first-order finite elements
and prove that the discrete energy functional converges to the continuous one
in the sense of $\Gamma$-convergence~\cite{braides2002}.
Moreover,
we propose two iterative algorithms
for the computation of low energy stationary points
based on time discretizations of the gradient flow of the energy functional
(see Algorithm~\ref{alg:coupled} and Algorithm~\ref{alg:decoupled} below).
These two algorithms differ from each other in the time discretization
(fully implicit for Algorithm~\ref{alg:coupled}, semi-implicit for Algorithm~\ref{alg:decoupled}).
For both algorithms, we prove well-posedness of the iteration,
an energy-decreasing property,
termination of the iterative loop,
an upper bound for the error in the unit-length constraint,
and
(under a restrictive assumption on the coefficients appearing in the energy functional)
convergence toward a stationary point.
Moreover, we perform numerical experiments to compare the two algorithms and to assess their performance.

Then, we extend the best performing algorithm (and its analysis) to the dynamic problem
and show that the resulting integrator (Algorithm~\ref{alg:tps}) is well-posed, stable,
and generates approximations that are unconditionally convergent toward a weak solution
of the coupled system of LLG equations.
A by-product of our constructive convergence proof
is the first proof of existence of weak solutions for this problem.

In general, the mathematical literature on AFM and FiM materials is much less developed
than that of FM materials.
We refer, e.g., to~\cite{bcko2021,bcko2022} for works 
discussing discrete-to-continuum variational limits of a two-dimensional atomistic model
of AFM materials.
As far as the continuum model considered in this work is concerned,
to the authors' knowledge, the only other work addressing it is~\cite{lcdw2020},
where extensions of the Gauss--Seidel projection method~\cite{wge2001,lxdcw2020}
have been proposed for its numerical approximation (but no convergence analysis is discussed).

To sum up, the novel contributions of the present work are the following:\\
$\bullet$
We provide the first mathematically rigorous formulation
of a state-of-the-art model currently used by applied scientists 
to simulate processes and devices involving AFM and FiM materials.\\
$\bullet$
Extending the techniques that have been developed for the approximation of the classical model of FM materials,
we introduce and analyze the first convergent numerical schemes for this model of AFM and FiM materials.

The remainder of the paper is organized as follows:
In Section~\ref{sec:model},
we present the mathematical model of AFM and FiM materials.
In Section~\ref{sec:preliminaries}, we introduce the main ingredient of our discretization.
The algorithms for the static problem, their properties, and two numerical experiments
are discussed in Section~\ref{sec:static}.
In Section~\ref{sec:llg}, we extend one of the algorithm and its analysis to the dynamic case,
and use it to simulate the dynamics of magnetic skyrmions in antiferromagnets.
In Section~\ref{sec:proofs}, we collect the proofs of all results of the work.
Finally, in Appendix~\ref{sec:nondimensional},
we show how to pass from the formulation of the model in physical units
to the dimensionless setting considered in this work.

\section{Mathematical model} \label{sec:model}

Let $\Omega\subset\R^3$ be a bounded Lipschitz domain
representing the volume occupied by an AFM or a FiM material.
The magnetic state of the material is described in terms of two unit-length vector fields,
$\mm_1$ and $\mm_2$.
The total magnetization of the sample is given by
$\mm = \eta_{\mathrm{s},1} \mm_1 + \eta_{\mathrm{s},2} \mm_2$,
where $\eta_{\mathrm{s},1},\eta_{\mathrm{s},2}>0$ are dimensionless constants.

In what follows,
for Lebesgue, Sobolev, and Bochner spaces and norms, we will use the standard notation~\cite{evans2010}.
To denote (spaces of) vector-valued or matrix-valued functions, we use bold letters, e.g.,
for any domain $U$, we denote both $L^2(U;\R^3)$ and $L^2(U;\R^{3 \times 3})$ by $\LL^2(U)$.
Moreover, we will denote by
$\inner{\cdot}{\cdot}$ the inner product of $\LL^2(\Omega)$
and by $\norm{\cdot}$ the corresponding norm
(any other inner product or norm will be denoted by the same notation,
but supplemented with a suitable subscript).
We will denote by $\inner{\cdot}{\cdot}$ also the duality pairing between $\HH^1(\Omega)$
and its dual, and note that it coincides with the inner product of $\LL^2(\Omega)$
if the arguments are in $\LL^2(\Omega)$.

\subsection{Static problem}
Stable magnetic configurations of the sample are described by minimizers
$\mm_1, \mm_2 : \Omega \to \sphere$
of the energy functional
\begin{equation} \label{eq:energy}
\E[\mm_1,\mm_2]
= \frac{1}{2} \sum_{\ell=1}^2 a_{\ell\ell} \norm{\Grad\mm_\ell}^2
+ a_{12} \inner{\Grad\mm_1}{\Grad\mm_2}
- a_0 \inner{\mm_1}{\mm_2},
\end{equation}
where the material constants $a_{11}, a_{22}, a_{12}, a_0 \in \R$
satisfy the inequalities
\begin{equation} \label{eq:assumption_coeff}
a_{11} + a_{22}>0 \quad \text{and} \quad a_{11}a_{22} > a_{12}^2.
\end{equation}
The three contributions in~\eqref{eq:energy} are called
\emph{inhomogeneous intralattice exchange},
\emph{inhomogeneous interlattice exchange},
and~\emph{homogeneous interlattice exchange},
respectively~\cite{spkcf2020}.
Minimizers are sought in the set of admissible pairs of vector fields
\begin{equation} \label{eq:admissible_set}
\begin{split}
\XXX
:&=
\HH^1(\Omega ; \sphere) \times \HH^1(\Omega ; \sphere) \\
&=
\{ (\mm_1 , \mm_2) \in \HH^1(\Omega) \times \HH^1(\Omega) : \abs{\mm_1} = \abs{\mm_2} = 1 \text{ a.e.\ in } \Omega \}.
\end{split}
\end{equation}
Note that~\eqref{eq:assumption_coeff} guarantees that
the energy is bounded from below in $\XXX$,
as there holds the inequality
$\E[\mm_1,\mm_2] \ge - \abs{a_0} \abs{\Omega}$
for all $(\mm_1, \mm_2) \in \XXX$,
and that the energy functional is weakly sequentially lower semicontinuous in $\HH^1(\Omega) \times \HH^1(\Omega)$
(see Proposition~\ref{prop:wlsc_energy} below).
Hence, existence of minimizers follows from the direct method of calculus of variations.

Stationary points of the energy are admissible pairs $(\mm_1, \mm_2) \in \XXX$
which, for all $\ell=1,2$, solve
\begin{equation} \label{eq:euler_lagrange1}
- \dual{\heff{\ell}[\mm_1,\mm_2]}{\vvphi - (\mm_\ell\cdot\vvphi)\mm_\ell} = 0
\quad
\text{for all }
\vvphi \in \HH^1(\Omega) \cap \LL^\infty(\Omega),
\end{equation}
where
the effective field $\heff{\ell}[\mm_1,\mm_2]$
is the (negative) Gateaux derivative of the energy with respect to $\mm_\ell$, i.e.,
\begin{equation} \label{eq:heff}
\begin{split}
\dual{\heff{\ell}[\mm_1,\mm_2]}{\pphi}
:&\stackrel{\phantom{\eqref{eq:energy}}}{=}
\bigg\langle - \frac{\delta\E[\mm_1,\mm_2]}{\delta \mm_{\ell}} , \pphi \bigg\rangle\\
&\stackrel{\eqref{eq:energy}}{=} - a_{\ell\ell} \inner{\Grad\mm_{\ell}}{\Grad\pphi}
- a_{12} \inner{\Grad\mm_{3-\ell}}{\Grad\pphi}
+ a_0 \inner{\mm_{3 - \ell}}{\pphi}.
\end{split}
\end{equation}
Equivalently, a stationary point $(\mm_1, \mm_2) \in \XXX$ can be seen as the solution of
\begin{equation} \label{eq:euler_lagrange2}
- \dual{\heff{\ell}[\mm_1,\mm_2]}{\pphi} = 0
\quad
\text{for all }
\pphi \in \KKK[\mm_\ell],
\end{equation}
where
\begin{equation} \label{eq:tangent_space}
\KKK[\mm_\ell] =
\{ \ppsi \in \HH^1(\Omega) : \mm_\ell \cdot \ppsi = 0 \text{ a.e.\ in } \Omega \}.
\end{equation}
Note that~\eqref{eq:euler_lagrange1} and~\eqref{eq:euler_lagrange2} can be interpreted as variational formulations
of the boundary value problem
\begin{equation} \label{eq:euler_lagrange_strong}
\begin{aligned}
- \mm_\ell \times \heff{\ell}[\mm_1,\mm_2] &= \0  \quad \text{in } \Omega,\\
\partial_{\nnu} \mm_\ell &= \0 \quad \text{on } \partial\Omega,\\
\end{aligned}
\end{equation}
where $\nnu: \partial\Omega \to \sphere$ denotes the outward-pointing unit normal vector to $\partial\Omega$.

\subsection{Dynamic problem}
Out of equilibrium,
the dynamics of the time-dependent vector fields $\mm_1, \mm_2: \Omega \times (0, \infty) \to \sphere$
is governed by a coupled system of two Landau--Lifshitz--Gilbert (LLG) equations,
one for each vector field:
\begin{equation} \label{eq:llg}
\mmt_\ell = - \eta_\ell \, \mm_\ell \times \heff{\ell}[\mm_1,\mm_2] + \alpha_\ell \, \mm_\ell \times \mmt_\ell
\quad
\text{for all } \ell = 1,2,
\end{equation}
where $\eta_\ell, \alpha_\ell>0$ are dimensionless constants.
Note that the two LLG equations are coupled to each other via their effective fields.
To complete the setting, \eqref{eq:llg} is supplemented with a suitable initial condition
and the same boundary conditions as in~\eqref{eq:euler_lagrange_strong}, i.e.,
\begin{equation} \label{eq:llg_ibc}
\mm_\ell(0) = \mm_\ell^0 \text{ in } \Omega
\quad\text{and}\quad
\partial_{\nnu}\mm_\ell = \0 \text{ on } \partial\Omega \times (0,\infty)
\quad
\text{for all } \ell = 1,2,
\end{equation}
for some admissible pair $(\mm_1^0, \mm_2^0) \in \XXX$.

In the following definition, we state the notion of a weak solution
to the initial boundary value problem~\eqref{eq:llg}--\eqref{eq:llg_ibc},
which naturally extends to the present setting the notion introduced in~\cite{as1992}
for the standard LLG equation.

\begin{definition}[weak solution] \label{def:weak}
Let $(\mm_1^0, \mm_2^0) \in \XXX$.
A global weak solution of~\eqref{eq:llg}--\eqref{eq:llg_ibc} is
$(\mm_1,\mm_2) \in L^{\infty}(0,\infty;\XXX)$
such that,
for all $T>0$, the following properties are satisfied:
\begin{itemize}
\item[\rm(i)] $\mm_{\ell}\vert_{\Omega_T} \in \HH^1(\Omega_T)$ for all $\ell=1,2$, where $\Omega_T := \Omega \times (0,T)$;
\item[\rm(ii)] $\mm_{\ell}(0)=\mm_{\ell}^0$ in the sense of traces for all $\ell=1,2$;
\item[\rm(iii)] for all $\ell=1,2$, for all $\vvphi\in\HH^1(\Omega_T)$,
it holds that
\begin{equation} \label{eq:weak:variational}
\begin{split}
& \int_0^T \inner{\mmt_{\ell}(t)}{\vvphi(t)} \, \dt \\
& \ = - \eta_\ell \int_0^T \dual{\heff{\ell}[\mm_1(t),\mm_2(t)]}{\vvphi(t) \times \mm_{\ell}(t)} \, \dt
+ \alpha_\ell \int_0^T \inner{\mm_{\ell}(t)\times\mmt_\ell(t)}{\vvphi(t)} \, \dt;
\end{split}
\end{equation}
\item[\rm(iv)] it holds that
\begin{equation} \label{eq:weak:energy}
\E[\mm_1(T),\mm_2(T)]
+ \sum_{\ell=1}^2 \frac{\alpha_\ell}{\eta_\ell} \int_0^T \norm{\mmt_{\ell}(t)}^2 \dt
\leq \E[\mm_1^0,\mm_2^0].
\end{equation}
\end{itemize}
\end{definition}

The variational formulations in~\eqref{eq:weak:variational}
are weak formulations of the LLG equations in~\eqref{eq:llg} in the space-time cylinder $\Omega_T$,
while~\eqref{eq:weak:energy} is a weak counterpart of the energy law
\begin{equation*}
\ddt\E[\mm_1(t),\mm_2(t)] = - \sum_{\ell=1}^2 \frac{\alpha_\ell}{\eta_\ell} \norm{\mmt_{\ell}(t)}^2 \leq 0
\quad
\text{for all } t > 0,
\end{equation*}
satisfied by sufficiently smooth solutions of~\eqref{eq:llg}.

\begin{remark}
For ease of presentation, we consider a dimensionless form of the energy functional.
We refer to Appendix~\ref{sec:physical} for its derivation
(starting from the equations in physical units usually encountered in the physical literature).
Moreover, we restrict ourselves to the case in which the energy comprises only the exchange contribution.
This is sufficient to capture the main mathematical features of the model:
First, the analytical and numerical treatment of standard lower-order energy contributions
(e.g., magnetocrystalline anisotropy, Zeeman energy, magnetostatic energy, Dzyaloshinskii--Moriya interaction)
is well understood (see, e.g., \cite{bffgpprs2014,hpprss2019}).
Second, lower-order terms do not entail the coupling of the fields
(see Appendix~\ref{sec:lower} for more details).
Hence, even in the presence of lower-order terms, the Euler--Lagrange equations~\eqref{eq:euler_lagrange1}
and the system of LLG equations~\eqref{eq:llg} are exchange-coupled only.
\end{remark}

\section{Preliminaries} \label{sec:preliminaries}

In this section,
we collect the notation and the definitions that are necessary to introduce our numerical schemes.

For the time discretization,
we consider uniform partitions of the positive real axis
with constant time-step size $\tau>0$, i.e., $t_i := i \tau$ for all $i \in \N_0$.
Given a sequence $\{ \phi^i\}_{i \in \N_0}$,
for all $i \in \N_0$ we define
$d_t \phi^{i+1} := (\phi^{i+1} - \phi^i)/\tau$.
We consider the time reconstructions $\phi_{\tau}$, $\phi^-_{\tau}$, $\phi^+_{\tau}$
defined, for all $i \in \N_0$ and $t \in [t_i,t_{i+1})$, as
\begin{equation} \label{eq:reconstructions}
\begin{split}
& \phi_{\tau}(t) := \frac{t-t_i}{\tau}\phi^{i+1} + \frac{t_{i+1} - t}{\tau}\phi^i,
\quad
\phi_{\tau}^-(t) := \phi^i,
\quad
\text{and}
\quad
\phi_{\tau}^+(t) := \phi^{i+1}.
\end{split}
\end{equation}
Note that $\partial_t \phi_\tau (t) = d_t \phi^{i+1}$ for all $i \in \N_0$ and $t \in [t_i,t_{i+1})$.

The spatial discretization is based on first-order finite elements.
We assume $\Omega$ to be a polyhedral domain
and consider a family $\{ \T_h \}_{h>0}$ of shape-regular tetrahedral meshes of $\Omega$ 
parametrized by the mesh size $h = \max_{K \in \T_h} \diam(K)$.
We denote by $\NN_h$ the set of vertices of $\T_h$.
For any $K \in \T_h$, let $\mathcal{P}^1(K)$ be the space of polynomials of degree at most $1$ on $K$.
We denote by $\S^1(\T_h)$ the space of piecewise affine and globally continuous functions
from $\Omega$ to $\R$, i.e.\
\begin{equation*}
\S^1(\T_h)
= \left\{v_h \in C^0(\overline{\Omega}): v_h \vert_K \in \mathcal{P}^1(K) \text{ for all } K \in \T_h \right\}.
\end{equation*}
It is well known that $\S^1(\T_h)$ is a finite-dimensional subspace of $H^1(\Omega)$ with $\dim\S^1(\T_h) = N_h := \#\NN_h$.
Let $\interp: C^0(\overline{\Omega}) \to \S^1(\T_h)$
denote the nodal interpolation operator,
i.e., for all $v \in C^0(\overline{\Omega})$, $\interp[v]$ is the unique element of $\S^1(\T_h)$ satisfying
$\interp[v](z) = v(z)$ for all $z \in \NN_h$.
We use the same notation to denote its vector-valued counterpart $\interp: \CC^0(\overline{\Omega}) \to \S^1(\T_h)^3$,
where the scalar-valued operator is applied to each component of a vector-valued function.

We consider the mass-lumped $L^2$-product $\inner[h]{\cdot}{\cdot}$ defined by
\begin{equation} \label{eq:mass-lumping}
\inner[h]{\ppsi}{\pphi}
= \int_\Omega \interp[\ppsi \cdot \pphi]
\quad \text{for all } \ppsi, \pphi \in \CC^0(\overline{\Omega}).
\end{equation}
We recall that this defines 
an inner product on $\S^1(\T_h)^3$
and that
the induced norm $\norm[h]{\cdot}$ satisfies the norm equivalence
\begin{equation} \label{eq:mass_lumping_equivalence}
\norm{\pphi_h}
\leq \norm[h]{\pphi_h}
\leq \sqrt{5} \, \norm{\pphi_h}
\quad \text{for all } \pphi_h \in \S^1(\T_h)^3,
\end{equation}
see \cite[Lemma~3.9]{bartels2015}.
Moreover, we have that
\begin{equation} \label{eq:h-inner-product}
\lvert \inner{\pphi_h}{\ppsi_h} -  \inner[h]{\pphi_h}{\ppsi_h} \rvert
\le C h^2 \norm[\LL^2(\Omega)]{\Grad\pphi_h} \norm[\LL^2(\Omega)]{\Grad\ppsi_h}
\quad \text{for all } \pphi_h, \ppsi_h \in \S^1(\T_h)^3,
\end{equation}
where $C>0$ depends only on the shape-regularity of $\T_h$;
see again~\cite[Lemma~3.9]{bartels2015}.

We conclude this section with a notational remark:
In what follows,
we will always denote by $C>0$ a generic constant, which will be always independent of the discretization parameters,
but not necessarily the same at each occurrence.

\section{Numerical energy minimization} \label{sec:static}

In this section,
we introduce a finite element discretization of the energy minimization problem
and show its convergence
in the sense of $\Gamma$-convergence.
Then,
we introduce two fully discrete algorithms to approximate stationary points of the energy
functional~\eqref{eq:energy}.
We state our results regarding well-posedness, stability, and convergence of the algorithms
and underpin our theoretical results with numerical experiments.
To make the presentation of the results concise,
all proofs are postponed to Section~\ref{sec:proofs_min}.

\subsection{Finite element discretization} \label{sec:gamma_convergence}

To discretize the set of admissible pairs in~\eqref{eq:admissible_set}, given a mesh $\T_h$ and a parameter $\delta>0$,
we consider the set
\begin{multline*}
\XXX_{h,\delta} := 
\{ (\mm_{1,h,\delta} , \mm_{2,h,\delta}) \in \S^1(\T_h)^3 \times \S^1(\T_h)^3 : \text{ for all } \ell=1,2,\\
\abs{\mm_{\ell,h,\delta}(z)} \ge 1 \text{ for all } z \in \NN_h
\text{ and }
\norm[L^1(\Omega)]{\interp[\abs{\mm_{\ell,h,\delta}}^2]-1} \le \delta \}.
\end{multline*}
Note that, at the discrete level, the unit-length constraint is relaxed~\cite{bartels2016,ahpprs2014}
and a mild control of the error is enforced by the inequality involving the parameter $\delta$.

The discrete static problem consists of seeking a minimizer
of the energy functional~\eqref{eq:energy} in the set of discrete admissible pairs in $\XXX_{h,\delta}$.
In the following theorem,
we show that the discrete energy functional $\E_{h,\delta}[\cdot,\cdot] := \E\vert_{\XXX_{h,\delta}}[\cdot,\cdot]$
converges toward the continuous one in the sense of $\Gamma$-convergence.
We note that our discretization is consistent, i.e., we do not modify the energy functional,
but we restrict the set in which minimizers are sought.

\begin{theorem}[$\Gamma$-convergence] \label{thm:gamma_convergence}
The following two properties hold:\\
{\rm(i)} Lim-inf inequality:
For every sequence $\{ (\mm_{1,h,\delta} , \mm_{2,h,\delta}) \}$
with $(\mm_{1,h,\delta} , \mm_{2,h,\delta}) \in \XXX_{h,\delta}$ for all $h,\delta>0$
such that,
for some $(\mm_1,\mm_2) \in \XXX$,
$\mm_{\ell,h,\delta} \weakto \mm_\ell$ in $\HH^1(\Omega)$ as $h,\delta \to 0$ for all $\ell =1,2$,
we have that
$\E[\mm_1,\mm_2] \le \liminf_{h,\delta \to 0} \E_{h,\delta}[\mm_{1,h,\delta} , \mm_{2,h,\delta}]$.\\
{\rm(ii)} Existence of a recovery sequence:
For every $(\mm_1,\mm_2) \in \XXX$,
there exists a sequence
$\{ (\mm_{1,h,\delta} , \mm_{2,h,\delta}) \}$
with $(\mm_{1,h,\delta} , \mm_{2,h,\delta}) \in \XXX_{h,\delta}$ for all $h,\delta>0$
such that
$(\mm_{1,h,\delta} , \mm_{2,h,\delta}) \to (\mm_1,\mm_2)$ in $\HH^1(\Omega) \times \HH^1(\Omega)$
and $\E_{h,\delta}[\mm_{1,h,\delta} , \mm_{2,h,\delta}] \to \E[\mm_1 , \mm_2]$
as $h,\delta \to 0$.
\end{theorem}

A well-known consequence of $\Gamma$-convergence
is the convergence of discrete global minimizers.

\begin{corollary} \label{cor:gamma_convergence}
Let $\{ (\mm_{1,h,\delta} , \mm_{2,h,\delta}) \}$ be a sequence 
such that $(\mm_{1,h,\delta} , \mm_{2,h,\delta}) \in \XXX_{h,\delta}$ is a global minimizer of
the discrete energy functional $\E_{h,\delta}[\cdot,\cdot]$ for all $h,\delta>0$.
Then, every accumulation point $(\mm_1,\mm_2)$ of the sequence
belongs to $\XXX$
and is a global minimizer of the continuous energy functional $\E[\cdot,\cdot]$.
\end{corollary}

We omit the proof of Corollary~\ref{cor:gamma_convergence}
as it is based on standard $\Gamma$-convergence arguments;
see, e.g., \cite[Section~1.5]{braides2002}.
Moreover, we recall that $\Gamma$-convergence
does not imply the convergence of local minimizers.

\subsection{Computation of low energy stationary points} \label{sec:results_static}

Let $\HHH$ be a Hilbert space with inner product $\inner[\HHH]{\cdot}{\cdot}$
such that $\XXX$ is continuously embedded in $\HHH \times \HHH$.
Furthermore, we suppose that there exists a constant $\cH \ge 1$ such that
\begin{equation} \label{eq:metric}
\cH^{-1} \norm{\pphi}
\le
\norm[\HHH]{\pphi}
\le
\cH \norm[\HH^1(\Omega)]{\pphi}
\quad
\text{for all } \pphi \in \HH^1(\Omega).
\end{equation}
To find stationary points with low energy,
we propose two iterative algorithms that are based
on two discretizations of
the dissipative dynamics governed by the $\HHH$-gradient flow of the energy
\begin{equation} \label{eq:gradient_flow}
\inner[\HHH]{\mmt_\ell}{\pphi}
+
\bigg\langle \frac{\delta\E[\mm_1,\mm_2]}{\delta \mm_{\ell}} , \pphi \bigg\rangle = \0
\quad
\text{for all } \pphi \in \KKK[\mm_\ell] \quad (\ell=1,2).
\end{equation}
The spatial discretization of both methods is based on first-order finite elements
as described in Section~\ref{sec:preliminaries}.
As a discrete counterpart of the space of pointwise orthogonal vector fields in~\eqref{eq:tangent_space},
for $\mm_h \in \S^1(\T_h)^3$ with $\mm_h(z) \neq \0$ for all $z \in \NN_h$, 
we consider the finite-dimensional space
\begin{equation} \label{eq:discrete_tangent_space}
\KKK_h[\mm_h]
:= \left\{
\pphi_h \in \S^1(\T_h)^3 : \mm_h(z) \cdot \pphi_h(z) = 0 \text{ for all } z \in \NN_h
\right\}.
\end{equation}
For discrete functions,
the pointwise orthogonality of~\eqref{eq:tangent_space} is required to hold only at the vertices of the mesh.
Note that $\KKK_h[\mm_h]$ is a subspace of $\S^1(\T_h)^3$ with dimension $2N_h$.
The time discretization is based on two different time-stepping methods.

\begin{remark}
In this section, we refer to the variable $t$ as \emph{time}
(accordingly, we refer to $\tau$ below as the \emph{time-step size}).
However, note that we are considering the static setting,
with the time variable $t$ playing the role of a \emph{pseudo-time}, introduced only for numerical purposes.
\end{remark}

The first method is proposed in the following algorithm.

\begin{algorithm}[coupled discrete gradient flow] \label{alg:coupled}
\textbf{Discretization parameters:}
Mesh size $h>0$, time-step size $\tau>0$, tolerance $\eps>0$.\\
\textbf{Input:}
Initial guess
$(\mm_{1,h}^0, \mm_{2,h}^0) \in \S^1(\T_h)^3 \times \S^1(\T_h)^3$ such that,
for all $\ell=1,2$,
$\abs{\mm_{\ell,h}^0(z)} = 1$ for all $z \in \NN_h$.\\
\textbf{Loop:}
For all $i \in \N_0$, iterate {\rm(i)--\rm(ii)} until the stopping criterion {\rm(stop)} is met:
\begin{itemize}
\item[\rm(i)] Given $(\mm_{1,h}^i, \mm_{2,h}^i) \in \S^1(\T_h)^3 \times \S^1(\T_h)^3$,
compute $(\vv_{1,h}^i, \vv_{2,h}^i) \in \KKK_h[\mm_{1,h}^i] \times \KKK_h[\mm_{2,h}^i]$ such that,
for all $(\pphi_{1,h} , \pphi_{2,h}) \in \KKK_h[\mm_{1,h}^i] \times \KKK_h[\mm_{2,h}^i]$
and $\ell=1,2$, it holds that
\begin{multline} \label{eq:coupled1}
\inner[\HHH]{\vv_{\ell,h}^i}{\pphi_{\ell,h}}
+ a_{\ell\ell} \tau \inner{\Grad\vv_{\ell,h}^i}{\Grad\pphi_{\ell,h}}
+ \frac{a_{12}}{2} \tau \inner{\Grad\vv_{3-\ell,h}^i}{\Grad\pphi_{\ell,h}}
- \frac{a_0}{2} \tau \inner{\vv_{3-\ell,h}^i}{\pphi_{\ell,h}} \\
=
- a_{\ell\ell} \inner{\Grad\mm_{\ell,h}^i}{\Grad\pphi_{\ell,h}}
- a_{12} \inner{\Grad\mm_{3-\ell,h}^i}{\Grad\pphi_{\ell,h}}
+ a_0 \inner{\mm_{3-\ell,h}^i}{\pphi_{\ell,h}}.
\end{multline}
\item[\rm(ii)] Define
\begin{equation} \label{eq:coupled2}
\mm_{\ell,h}^{i+1} := \mm_{\ell,h}^i + \tau \vv_{\ell,h}^i
\quad \text{for all } \ell =1,2.
\end{equation}
\item[\rm(stop)]
Stop iterating {\rm(i)--(ii)} if $(\vv_{1,h}^i, \vv_{2,h}^i) \in \KKK_h[\mm_{1,h}^i] \times \KKK_h[\mm_{2,h}^i]$ satisfies
\begin{equation} \label{eq:coupled_stopping}
\max_{\ell=1,2} \left( \norm[\HHH]{\vv_{\ell,h}^i}^2 + \tau \norm{\Grad\vv_{\ell,h}^i}^2 \right)
\le \eps^2 \abs{\Omega}.
\end{equation}
\end{itemize}
\textbf{Output:}
If $i^* \in \N_0$ denotes the smallest integer satisfying the stopping criterion~\eqref{eq:coupled_stopping},
define the approximate stationary point
$(\mm_{1,h},\mm_{2,h}) := (\mm_{1,h}^{i^*},\mm_{2,h}^{i^*})$.
\end{algorithm}

The second method is proposed in the following algorithm.

\begin{algorithm}[decoupled discrete gradient flow] \label{alg:decoupled}
\textbf{Discretization parameters:} Mesh size $h>0$,
time-step size $\tau>0$, tolerance $\eps>0$.\\
\textbf{Input:}
Initial guess
$(\mm_{1,h}^0, \mm_{2,h}^0) \in \S^1(\T_h)^3 \times \S^1(\T_h)^3$ such that,
for all $\ell=1,2$,
$\abs{\mm_{\ell,h}^0(z)} = 1$ for all $z \in \NN_h$.\\
\textbf{Loop:}
For all $i \in \N_0$, iterate {\rm(i)--\rm(ii)} until the stopping criterion {\rm(stop)} is met:
\begin{itemize}
\item[\rm(i)] Given $(\mm_{1,h}^i, \mm_{2,h}^i) \in \S^1(\T_h)^3 \times \S^1(\T_h)^3$,
compute $(\vv_{1,h}^i, \vv_{2,h}^i) \in \KKK_h[\mm_{1,h}^i] \times \KKK_h[\mm_{2,h}^i]$ such that,
for all $(\pphi_{1,h} , \pphi_{2,h}) \in \KKK_h[\mm_{1,h}^i] \times \KKK_h[\mm_{2,h}^i]$
and $\ell=1,2$, it holds that
\begin{multline} \label{eq:decoupled1}
\inner[\HHH]{\vv_{\ell,h}^i}{\pphi_{\ell,h}}
+ a_{\ell\ell} \tau \inner{\Grad\vv_{\ell,h}^i}{\Grad\pphi_{\ell,h}} \\
=
- a_{\ell\ell} \inner{\Grad\mm_{\ell,h}^i}{\Grad\pphi_{\ell,h}}
- a_{12} \inner{\Grad\mm_{3-\ell,h}^i}{\Grad\pphi_{\ell,h}}
+ a_0 \inner{\mm_{3-\ell,h}^i}{\pphi_{\ell,h}}.
\end{multline}
\item[\rm(ii)] Define
\begin{equation} \label{eq:decoupled2}
\mm_{\ell,h}^{i+1} := \mm_{\ell,h}^i + \tau \vv_{\ell,h}^i
\quad \text{for all } \ell =1,2.
\end{equation}
\item[\rm(stop)]
Stop iterating {\rm(i)--(ii)} if $(\vv_{1,h}^i, \vv_{2,h}^i) \in \KKK_h[\mm_{1,h}^i] \times \KKK_h[\mm_{2,h}^i]$ satisfies
\begin{equation} \label{eq:decoupled_stopping}
\max_{\ell=1,2} \left( \norm[\HHH]{\vv_{\ell,h}^i}^2 + \tau \norm{\Grad\vv_{\ell,h}^i}^2 \right)
\le \eps^2 \abs{\Omega}.
\end{equation}
\end{itemize}
\textbf{Output:}
If $i^* \in \N_0$ denotes the smallest integer satisfying the stopping criterion~\eqref{eq:decoupled_stopping},
define the approximate stationary point
$(\mm_{1,h},\mm_{2,h}) := (\mm_{1,h}^{i^*},\mm_{2,h}^{i^*})$.
\end{algorithm}

In both Algorithm~\ref{alg:coupled} and Algorithm~\ref{alg:decoupled} the iteration stops
when the size of the updates is sufficiently small
(according to \eqref{eq:coupled_stopping} and \eqref{eq:decoupled_stopping}, respectively).
The algorithms
are characterized by a different treatment of the
inhomogeneous and homogeneous interlattice exchange contributions,
which are treated implicitly in Algorithm~\ref{alg:coupled}
and explicitly in Algorithm~\ref{alg:decoupled}.
One immediate consequence is that in Algorithm~\ref{alg:coupled} the two equations are coupled
(as they are in the continuous problem) and
one iteration of the algorithm requires the solution of \emph{one} $4N_h$-by-$4N_h$ linear system,
whereas in Algorithm~\ref{alg:decoupled} the two equations are decoupled
and one iteration of the algorithm requires the solution of \emph{two} $2N_h$-by-$2N_h$ linear systems
(that are independent of each other and thus can be solved in parallel).
This difference will affect the solvability and energetic behavior of the algorithms,
which will be the subject of the following propositions.

In the following proposition,
we establish the properties of Algorithm~\ref{alg:coupled}.

\begin{proposition}[properties of Algorithm~\ref{alg:coupled}] \label{prop:coupled}
There hold the following statements:\\
{\rm(i)}
Suppose that $\tau$ satisfies $\cH^2 \abs{a_0} \tau < 2$,
where $a_0$ is one of the coefficients in~\eqref{eq:energy}
and $\cH$ is the constant in~\eqref{eq:metric}.
Then, for all $i \in \N_0$, \eqref{eq:coupled1} admits a unique solution
$(\vv_{1,h}^i, \vv_{2,h}^i) \in \KKK_h[\mm_{1,h}^i] \times \KKK_h[\mm_{2,h}^i]$.\\
{\rm(ii)}
Under the assumption of part~{\rm(i)},
suppose that $\tau$ additionally satisfies $\cH \cT \abs{a_0} \tau < 1$,
where
$\cT > 0$ is a constant which depends only on the shape-regularity of the family of meshes.
Then, Algorithm~\ref{alg:coupled} terminates within a finite number of iterations.
In particular, the approximate stationary point $(\mm_{1,h},\mm_{2,h})$ is well defined.\\
{\rm(iii)}
Under the assumption of part~{\rm(i)},
for all $i \in \N_0$, the iterates of Algorithm~\ref{alg:coupled} satisfy
\begin{equation} \label{eq:coupled_energy_law}
\E[\mm_{1,h}^{i+1},\mm_{2,h}^{i+1}]
- \E[\mm_{1,h}^i,\mm_{2,h}^i]
=
- \tau \sum_{\ell=1}^2 \norm[\HHH]{\vv_{\ell,h}^i}^2
- \frac{\tau^2}{2} \sum_{\ell=1}^2 a_{\ell\ell} \norm{\Grad\vv_{\ell,h}^i}^2.
\end{equation}
In particular,
the sequence of energies generated by the algorithm is monotonically decreasing,
i.e.,
it holds that $\E[\mm_{1,h}^{i+1},\mm_{2,h}^{i+1}] \le \E[\mm_{1,h}^i,\mm_{2,h}^i]$.\\
{\rm(iv)}
Under the assumptions of part~{\rm(ii)},
there exists $C>0$ such that
the approximate stationary point $(\mm_{1,h},\mm_{2,h})$ satisfies
\begin{equation*}
\norm[L^1(\Omega)]{\interp[\abs{\mm_{\ell,h}}^2]-1}
\le C \tau \left( 1 + \sum_{\ell=1}^2 \norm[\HH^1(\Omega)]{\mm_{\ell,h}^0}^2 \right)
\quad
\text{for all } \ell=1,2.
\end{equation*}
The constant $C>0$ depends only on $a_{11}$, $a_{12}$, $a_{22}$, $a_0$, $\cH$,
and the shape-regularity of the family of meshes.
\end{proposition}

Corresponding results for Algorithm~\ref{alg:decoupled} are the subject of the following proposition.

\begin{proposition}[properties of Algorithm~\ref{alg:decoupled}] \label{prop:decoupled}
There hold the following statements:\\
{\rm(i)}
For all $i \in \N_0$,
\eqref{eq:decoupled1} admits a unique solution
$(\vv_{1,h}^i, \vv_{2,h}^i) \in \KKK_h[\mm_{1,h}^i] \times \KKK_h[\mm_{2,h}^i]$.\\
{\rm(ii)}
Suppose that $\tau$ satisfies $\cH (\cH/2 + \cT) \abs{a_0} \tau < 1$,
where $a_0$ is one of the coefficients in~\eqref{eq:energy},
$\cH$ is the constant in~\eqref{eq:metric},
and $\cT > 0$ is a constant which depends only on the shape-regularity of the family of meshes.
Then, Algorithm~\ref{alg:decoupled} terminates within a finite number of iterations.
In particular, the approximate stationary point $(\mm_{1,h},\mm_{2,h})$ is well defined.\\
{\rm(iii)}
For all $i \in \N_0$,
the iterates of Algorithm~\ref{alg:decoupled} satisfy
\begin{multline} \label{eq:decoupled_energy_law}
\E[\mm_{1,h}^{i+1},\mm_{2,h}^{i+1}]
- \E[\mm_{1,h}^i,\mm_{2,h}^i]
=
- \tau \sum_{\ell=1}^2 \norm[\HHH]{\vv_{\ell,h}^i}^2
- \frac{\tau^2}{2} \sum_{\ell=1}^2 a_{\ell\ell} \norm{\Grad\vv_{\ell,h}^{i}}^2 \\
+ a_{12} \tau^2 \inner{\Grad\vv_{1,h}^{i}}{\Grad\vv_{2,h}^{i}}
- a_0 \tau^2 \inner{\vv_{1,h}^{i}}{\vv_{2,h}^{i}}.
\end{multline}
Moreover,
if $\tau$ satisfies $\cH^2 \abs{a_0} \tau \le 2$,
then
the sequence of energies generated by the algorithm is monotonically decreasing,
i.e.,
it holds that $\E[\mm_{1,h}^{i+1},\mm_{2,h}^{i+1}] \le \E[\mm_{1,h}^i,\mm_{2,h}^i]$.\\
{\rm(iv)}
Under the assumptions of part~{\rm(ii)},
there exists $C>0$ such that
the approximate stationary point $(\mm_{1,h},\mm_{2,h})$ satisfies
\begin{equation*}
\norm[L^1(\Omega)]{\interp[\abs{\mm_{\ell,h}}^2]-1}
\le C \tau \left( 1 + \sum_{\ell=1}^2 \norm[\HH^1(\Omega)]{\mm_{\ell,h}^0}^2 \right)
\quad
\text{for all } \ell=1,2.
\end{equation*}
The constant $C>0$ depends only on $a_{11}$, $a_{12}$, $a_{22}$, $a_0$, $\cH$,
and the shape-regularity of the family of meshes.
\end{proposition}

Each iteration of Algorithm~\ref{alg:coupled} is well defined if the time-step size is sufficiently small.
Moreover, the algorithm unconditionally generates a monotonically decreasing sequence of energies.
Conversely, each iteration of Algorithm~\ref{alg:decoupled} is unconditionally well defined,
but the sequence of energies it generates is monotonically decreasing only if the time-step size is sufficiently small.
Furthermore, we note that the inequalities in point~(iv)
of both Proposition~\ref{prop:coupled} and Proposition~\ref{prop:decoupled}
show that if the initial guesses are uniformly bounded in $\HH^1(\Omega)$
(in the sense of~\eqref{eq:uniform_boundedness_initial} below),
then the approximate stationary points generated by the algorithms
belong to the set of admissible pairs $\XXX_{h,\delta}$ with $\delta$ of the form $\delta = C \tau$.

In the following theorem,
we show that
the sequence of approximate stationary points computed with both algorithms
converges toward an admissible pair in $\XXX$ as the discretization parameters go to zero.
If we neglect the inhomogeneous interlattice exchange contribution,
we can identify the limit with a stationary point of the energy functional~\eqref{eq:energy}.

\begin{theorem}[convergence of Algorithm~\ref{alg:coupled} and Algorithm~\ref{alg:decoupled}] \label{thm:static}
Suppose that there exists $c_0>0$, independent of the discretization parameters $h$, $\tau$, and $\eps$,
such that
\begin{equation} \label{eq:uniform_boundedness_initial}
\sup_{h>0} \left( \sum_{\ell=1}^2 \norm[\HH^1(\Omega)]{\mm_{\ell,h}^0}^2 \right) \le c_0.
\end{equation}
Suppose that $\tau \to 0$ and $\eps \to 0$ as $h \to 0$.
Then, as $h \to 0$, the sequence of
approximate stationary points $\{ (\mm_{1,h},\mm_{2,h}) \}_{h>0}$
generated by either Algorithm~\ref{alg:coupled} or Algorithm~\ref{alg:decoupled},
upon extraction of a subsequence,
converges weakly in $\HH^1(\Omega) \times \HH^1(\Omega)$
toward a point $(\mm_{1},\mm_{2}) \in \XXX$.
If $a_{12}=0$, the limit $(\mm_{1},\mm_{2}) \in \XXX$  is stationary point of the energy functional~\eqref{eq:energy}.
\end{theorem}

A byproduct of Theorem~\ref{thm:static} is the existence of weak solutions
to the Euler--Lagrange equations~\eqref{eq:euler_lagrange1} for the case $a_{12}=0$.

\begin{remark} \label{rem:why_assumption}
In our analysis, we can identify the limit of the sequence of approximate stationary points
with a stationary point of the energy only if we assume that $a_{12}=0$,
i.e., if we neglect the inhomogeneous interlattice exchange contribution from the energy.
This restriction is related to the fact that,
if $a_{12} \neq 0$, the weak formulation of the approximate Euler--Lagrange equations
satisfied by $(\mm_{1,h},\mm_{2,h})$ contains a term that involves the $L^2$-product
of $\Grad\mm_{1,h}$ and $\Grad\mm_{2,h}$.
Since $(\mm_{1,h},\mm_{2,h})$ converges to $(\mm_1,\mm_2)$ only weakly in $\HH^1(\Omega) \times \HH^1(\Omega)$,
we are not allowed to pass this term to the limit.
We believe that this issue comes from the fact that our algorithms do not use any regularization,
so that the stability analysis does not yield any additional regularity
(and thus no stronger convergence properties)
that would allow us
to use arguments based on compensated compactness (see, e.g., \cite[Chapter~5]{evans1990} or \cite[Section~I.3]{struwe2008}).
However, we note that our numerical experiments suggest that the algorithms behave well
even if $a_{12} \neq 0$.
Moreover, in many situations (see, e.g., \cite{pkcatsf2019,sztfm2020}), the
inhomogeneous interlattice exchange contribution is of limited physical value
and is omitted, so that the current theory already covers many applications.
\end{remark}

\subsection{Numerical experiments} \label{sec:numerics_stationary}


Before moving to the dynamic case,
we show the effectivity of the proposed algorithms with two numerical experiments.
The computations presented in this section (and in Section~\ref{sec:numerics_llg} below) have been performed
with an implementation based on the open-source finite element library Netgen/NGSolve~\cite{netgen} (version 6.2.2302).
Lower-order energy contributions such as magnetocrystalline anisotropy, Dzyaloshinskii--Moriya interaction, and Zeeman energy (cf.\ Section~\ref{sec:lower}),
omitted in our analysis, are treated explicitly (and thus contribute only to the right-hand-sides of~\eqref{eq:coupled1} and~\eqref{eq:decoupled1}); see~\cite{bffgpprs2014,ahpprs2014}.
The orthogonality constraint in~\eqref{eq:coupled1} and~\eqref{eq:decoupled1} is enforced using the null-space method
discussed in~\cite{ramage2013preconditioned,kraus2019iterative}.
The resulting linear systems are solved using the generalized minimal residual method (GMRES)
with an incomplete LU decomposition preconditioner.
We note that in the static case, the use of the conjugate gradient method is possible due to symmetry,
but we use GMRES in these tests to maintain consistency with the dynamic case
(see Section~\ref{sec:numerics_llg} below).
All computations have been made on an i5-9500 CPU with \SI{16}{\giga\byte} of installed memory. 
Magnetization configurations are visualized with ParaView~\cite{agl2005}.

\subsubsection{Comparison of the algorithms} \label{sec:num_compare}

In this experiment, we aim to compare to each other
Algorithm~\ref{alg:coupled} and Algorithm~\ref{alg:decoupled},
and to evaluate the impact on their performance of the choice of the gradient flow metric,
i.e., the inner product $\inner[\HHH]{\cdot}{\cdot}$ in~\eqref{eq:gradient_flow}.

For the dimensionless setting discussed in Section~\ref{sec:model},
we consider a toy problem on the unit cube $\Omega = (-1/2,1/2)^3$.
The total energy consists of exchange and uniaxial anisotropy, i.e.,
\begin{multline*}
\E[\mm_1,\mm_2]
= \frac{1}{2} \sum_{\ell=1}^2 a_{\ell\ell} \int_\Omega \abs{\Grad\mm_\ell}^2
+ a_{12} \int_\Omega \Grad\mm_1 : \Grad\mm_2
- a_0 \int_\Omega \mm_1 \cdot \mm_2 \\
+ \frac{q_1^2}{2} \int_{\Omega} [1 - (\aa\cdot\mm_1)^2]
+ \frac{q_2^2}{2} \int_{\Omega} [1 - (\aa\cdot\mm_2)^2],
\end{multline*}
with exchange constants $a_{11} = 2$, $a_{22} = 1$, $a_{12} = -1/2$, and $a_0 = - 100$,
anisotropy constants $q_1 = 5$ and $q_2 = 10$,
and easy axis $\aa = (1,1,1)/\sqrt{3}$.
It is easy to see that for this setup the energy minimization problem admits two global minimizers
$(\mm_1^\pm,\mm_2^\pm) \equiv \pm (\aa, -\aa)$
and that the energy value at the minimizers is $\E[\mm_1^\pm,\mm_2^\pm]= -100$.

For the discretization,
we consider a tetrahedral mesh generated by Netgen with mesh size $h \approx 0.209$ (\num{1433} vertices and \num{6201} elements),
and we set $\tau = 10^{-3}$ and $\eps = 10^{-4}$.
Starting from the constant initial guess $\mm_{1,h}^0 \equiv (1,0,0)$ and $\mm_{2,h}^0 \equiv (0,1,0)$,
we run Algorithm~\ref{alg:coupled} and Algorithm~\ref{alg:decoupled}
for three different choices for the gradient flow metric:
the $L^2$-metric $\inner[\HHH]{\cdot}{\cdot} = \inner{\cdot}{\cdot}$,
the mass-lumped $L^2$-metric $\inner[\HHH]{\cdot}{\cdot} = \inner[h]{\cdot}{\cdot}$ (see~\eqref{eq:mass-lumping}),
and the $H^1$-metric $\inner[\HHH]{\cdot}{\cdot} = \inner{\cdot}{\cdot} + \inner{\Grad(\cdot)}{\Grad(\cdot)}$.
For all six runs (two algorithms, three metrics each), the iterative algorithm returns as approximate stationary point
an approximation of the minimizer $(\mm_1^-,\mm_2^-) \equiv (-\aa, \aa)$.

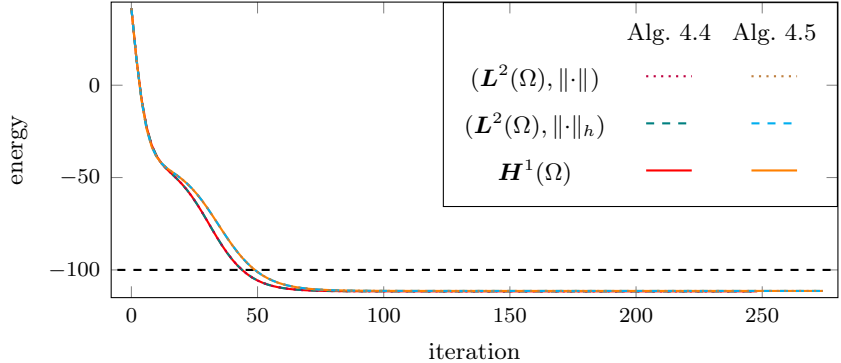
\begin{figure}[ht]
\centering
\begin{tikzpicture}
\begin{axis}
[
width = 0.7\textwidth,
height = 5.5cm,
xlabel = {\scriptsize iteration},
ylabel = {\scriptsize energy},
xmin = -8,
xmax = 280,
ymin = -115,
ymax = 45
]
\addplot [domain = -50:300,thick,dashed]{-100};
\addplot[red, thick]	table[x=t, y=energy, col sep=comma]{DataFiles/Constant_Coupled_H1.dat};
\label{plot:Cons_Coup_H1}
\addplot[teal, thick, dashed]	table[x=t, y=energy, col sep=comma]{DataFiles/Constant_Coupled_L2Mass.dat};
\label{plot:Cons_Coup_L2Mass}
\addplot[purple, thick, dotted]	table[x=t, y=energy, col sep=comma]{DataFiles/Constant_Coupled_L2.dat};
\label{plot:Cons_Coup_L2}
\addplot[orange, thick]	table[x=t, y=energy, col sep=comma]{DataFiles/Constant_Decoupled_H1.dat};
\label{plot:Cons_Decoup_H1}
\addplot[cyan, thick, dashed]	table[x=t, y=energy, col sep=comma]{DataFiles/Constant_Decoupled_L2Mass.dat};
\label{plot:Cons_Decoup_L2Mass}
\addplot[brown, thick, dotted]	table[x=t, y=energy, col sep=comma]{DataFiles/Constant_Decoupled_L2.dat};
\label{plot:Cons_Decoup_L2}
\coordinate (legend) at (axis description cs:1,1);
\end{axis}
\matrix [
draw,
matrix of nodes,
anchor=north east,
] at (legend) {
									& \scriptsize Alg.~\ref{alg:coupled}	& \scriptsize Alg.~\ref{alg:decoupled}\\
\scriptsize $(\LL^2(\Omega),\norm{\cdot})$	& \ref*{plot:Cons_Coup_L2}		& \ref*{plot:Cons_Decoup_L2}\\
\scriptsize $(\LL^2(\Omega),\norm[h]{\cdot})$	& \ref*{plot:Cons_Coup_L2Mass}	& \ref*{plot:Cons_Decoup_L2Mass}\\
\scriptsize $\HH^1(\Omega)$				& \ref*{plot:Cons_Coup_H1}		& \ref*{plot:Cons_Decoup_H1}\\
};
\end{tikzpicture}
\caption{Experiment of Section~\ref{sec:num_compare}:
Evolution of the energy with different algorithms and gradient flow metrics (constant initial guess).}
\label{fig:constant_energy_graph}
\end{figure}

In Table~\ref{tab:constant_initial_test},
we compare the performance of each combination in terms of
\begin{itemize}
\item
the final energy $\E[\mm_{1,h},\mm_{2,h}]$ of the approximate stationary point (\emph{energy});
\item
the difference $\E[\hat\mm_{1,h},\hat\mm_{2,h}]+100$ between the final energy $\E[\hat\mm_{1,h},\hat\mm_{2,h}]$ of the normalized approximate stationary point (\emph{proj.\ energy err.}) and the expected energy \num{-100},
where $\hat\mm_{\ell,h} = \interp[\mm_{\ell,h}/\abs{\mm_{\ell,h}}]$ for all $\ell=1,2$;
\item
the number of iterations necessary to meet the stopping criterion (\emph{num.\ iter.});
\item
the average solve time per iteration (\emph{solve time}), measured in \si{\second},
where the solve time is defined as the time needed to solve the linear system~\eqref{eq:coupled1} for Algorithm~\ref{alg:coupled}
and as the sum of the times needed to solve the two linear systems~\eqref{eq:decoupled1} (one for $\ell=1$ and one for $\ell=2$) for Algorithm~\ref{alg:decoupled};
\item
the error in the unit-length constraint measured in the $L^1$-norm, i.e.,
$\mathrm{err}_{L^1} := \max_{\ell=1,2} \big\lVert \interp\big[\abs{\mm_{\ell,h}}^2 \big]-1 \rVert_{L^1(\Omega)}$;
\item
the error in the unit-length constraint measured in the $L^\infty$-norm, i.e.,
$\mathrm{err}_{L^\infty} := \max_{\ell=1,2} \lVert \mm_{\ell,h} \rVert_{\LL^{\infty}(\Omega)} - 1$.
\end{itemize}
Moreover, in Figure~\ref{fig:constant_energy_graph}, for all six combinations of algorithms and metrics,
we plot the evolution of the energy during the iteration.

\begin{table}[ht]
\scriptsize
\begin{tabular}{c|c|c|c|c|c|c} 
& \multicolumn{3}{c|}{Algorithm~\ref{alg:coupled} (coupled)} & \multicolumn{3}{c}{Algorithm~\ref{alg:decoupled} (decoupled)}\\
\hline
$\HHH$  & $(\LL^2(\Omega),\norm{\cdot})$ & $(\LL^2(\Omega),\norm[h]{\cdot})$ & $\HH^1(\Omega)$ & $(\LL^2(\Omega),\norm{\cdot})$ & $(\LL^2(\Omega),\norm[h]{\cdot})$ & $\HH^1(\Omega)$\\ 
\hline
\emph{energy} & \num{-111.59} & \num{-111.59} & \num{-111.59} & \num{-111.38} & \num{-111.38} & \num{-111.38} \\
\emph{proj.\ energy err.} & \num{8.56e-11} & \num{8.56e-11} & \num{8.56e-11} & \num{9.90e-11} & \num{9.90e-11} & \num{9.90e-11} \\
\emph{num.\ iter.} & 249 & 249 & 249 & 275 & 275 & 275 \\
\emph{solve time} & 0.223 & 0.232 & 0.376 & 0.095 & 0.094 & 0.209 \\
$\mathrm{err}_{L^\infty}$ & 0.049 & 0.049 & 0.049 & 0.047 & 0.047 & 0.047 \\
$\mathrm{err}_{L^1}$ & \num{0.100} & \num{0.100} & \num{0.100} & \num{9.61e-02} & \num{9.61e-02} & \num{9.61e-02}
\end{tabular}
\caption{Experiment of Section~\ref{sec:num_compare}: Comparison of algorithms and gradient flow metrics (constant initial guess).\label{tab:constant_initial_test}}
\end{table}

In the very first part of the gradient flow dynamics (corresponding roughly to the first 15 iterations),
the constant initial guess with $\mm_1$ and $\mm_2$ perpendicular to each other evolves to reach a constant state with an antiparallel alignment of the fields.
This yields a strong reduction of the
$a_0$-modulated homogeneous interlattice exchange contribution
(with the total energy abruptly dropping from an initial value of about 41 to -42).
The rest of the dynamics is slower and consists in a rotation of the pair of constant fields
which make them align to the direction of the easy axis
as prescribed by the anisotropy energy contribution.

Looking at the results,
we observe that Algorithm~\ref{alg:coupled} and Algorithm~\ref{alg:decoupled} require approximately the same number of iterations
to fulfill the stopping criterion (those of Algorithm~\ref{alg:coupled} are slightly less than those of Algorithm~\ref{alg:decoupled}).
The energy decay of Algorithm~\ref{alg:coupled} is faster than the one of Algorithm~\ref{alg:decoupled},
but this does not lead to a significantly smaller number of iterations.
On the other hand, the average solve time of Algorithm~\ref{alg:decoupled} is about half of the one of Algorithm~\ref{alg:coupled},
which makes the simulations performed with the decoupled algorithm significantly faster.
The different metrics are practically identical
(except for a minimal difference in the average solve time).
For the $L^2$- and mass-lumped $L^2$-metric, this is unsurprising,
since they are equivalent to each other; see~\eqref{eq:mass_lumping_equivalence}.
We believe that the equivalence to the $H^1$-metric in this example (with constant initial guess) is due
to the fact that the updates $\vv_{\ell,h}^{i}$ are essentially uniform,
and hence their gradients  $\Grad \vv_{\ell,h}^{i}$ are essentially zero.
It follows that in the numerical scheme, the gradient part of the $H^1$-metric is small,
reducing to the $L^2$-metric.

There is a significant discrepancy between the value of the energy at the minimizer ($\E[\mm_{1}^+,\mm_{2}^+]=$ \num{-100})
and the one of its approximation ($\E[\mm_{1,h},\mm_{2,h}] \approx$ -111).
However, if we remove the error in the unit-length constraint by normalizing the fields at the vertices of the mesh,
we obtain the desired value up to the tenth digit.
This shows that our projection-free algorithms are perfectly able to identify the minimizers.
However, for a quantitative match of the energy values, the error in the constraint needs to be removed or reduced
(applying a nodal projection to the final configuration or decreasing the time-step size).

\begin{table}[ht]
\scriptsize
\begin{tabular}{c|c|c|c|c|c|c} 
& \multicolumn{3}{c|}{Algorithm~\ref{alg:coupled} (coupled)} & \multicolumn{3}{c}{Algorithm~\ref{alg:decoupled} (decoupled)}\\
\hline
$\HHH$  & $(\LL^2(\Omega),\norm{\cdot})$ & $(\LL^2(\Omega),\norm[h]{\cdot})$ & $\HH^1(\Omega)$ & $(\LL^2(\Omega),\norm{\cdot})$ & $(\LL^2(\Omega),\norm[h]{\cdot})$ & $\HH^1(\Omega)$\\ 
\hline
\emph{energy} & -182.36 & -158.04 & -100.38 & -182.70 & -158.03 & -100.38 \\
\emph{proj.\ energy err.} & \num{5.70e-11} & \num{6.33e-11} & \num{6.02e-9} & \num{6.54e-11} & \num{7.70e-11} & \num{6.02e-9} \\
\emph{num.\ iter.} & 231 & 231 & 13713 & 252 & 254 & 13718 \\
\emph{solve time} & 0.220 & 0.227 & 0.347 & 0.093 & 0.092 & 0.170 \\
$\mathrm{err}_{L^\infty}$ & 1.185 & 0.923 & 0.004 & 1.196 & 0.905 & 0.004\\
$\mathrm{err}_{L^1}$ & \num{0.668} & \num{0.433} & \num{2.51e-03} & \num{0.670} & \num{0.428} & \num{2.51e-03}\\
\end{tabular}
\caption{Experiment of Section~\ref{sec:num_compare}:
Comparison of algorithms and gradient flow metrics (random initial guess).}
\label{tab:random_initial_test}
\end{table}

Next, we repeat the experiment, but this time we start from a random initial guess (the same for all simulations).
For all six runs, the iterative algorithm again returns as approximate stationary point
an approximation of the minimizer $(\mm_1^-,\mm_2^-) \equiv (-\aa, \aa)$.
The results are displayed in Table~\ref{tab:random_initial_test}.
The faster average solve time of Algorithm~\ref{alg:decoupled} observed for the case of a constant initial guess is confirmed.
However, in this case,
we observe a clear difference between the $L^2$-metrics (with and without mass lumping) and the $H^1$-metric,
with the latter requiring a significantly larger number of iterations (ca 13700 against 200--300),
which results in longer computational times.
However, as far as the unit-length constraint is concerned, the $H^1$-metric is characterized by a much better accuracy.

Overall, our experiments show that the decoupled approach of Algorithm~\ref{alg:decoupled},
due to its computational efficiency, is preferable over the coupled one of Algorithm~\ref{alg:coupled}.
On the other hand, the choice of the gradient flow metric is more delicate.
While for a constant initial guess (with low exchange energy) the metrics are essentially equivalent,
for a random initial guess (with large exchange energy) the $H^1$-metric guarantees a significantly smaller violation
of the unit-length constraint at the discrete level (which, however, is obtained at the price of higher computational costs).

\subsubsection{Skyrmion formation} \label{sec:numerics_skyrmion}

In this experiment, we aim to highlight the capability of our algorithms
to compute stable magnetization configurations in AFM materials.

The domain is an AFM nanodisk of thickness \SI{1}{\nano\meter}
(aligned with $x_3$-axis) and diameter \SI{60}{\nano\meter} (aligned with the $x_1x_2$-plane).
The energy consists of exchange, out-of-plane uniaxial anisotropy, and interfacial Dzyaloshinskii--Moriya interaction,
and reads as
\begin{multline*}
\E[\mm_1,\mm_2]
= \frac{1}{2} \sum_{\ell=1}^2 a_{\ell\ell} \int_\Omega \abs{\Grad\mm_\ell}^2
- a_0 \int_\Omega \mm_1 \cdot \mm_2
+ \frac{q^2}{2} \int_{\Omega} [1 - (\aa\cdot\mm_1)^2] \\
+ \frac{q^2}{2} \int_{\Omega} [1 - (\aa\cdot\mm_2)^2]
+ \int_{\Omega} \widehat\DD : (\Grad\mm_1 \times \mm_1)
+ \int_{\Omega} \widehat\DD : (\Grad\mm_2 \times \mm_2),
\end{multline*}
where the dimensionless parameters $a_{11}$, $a_{22}$, $a_0$, $q$ and $\widehat\DD$
are obtained from the material parameters collected in Table~\ref{tab:material_parameters}
as explained in Appendix~\ref{sec:nondimensional}.

\begin{table}[ht]
\scriptsize
\setlength{\tabcolsep}{10pt}
\begin{tabular}{c|c}
Parameter & Value \\
\hline
$\Msl{1}$, $\Msl{2}$ & \SI{376}{\kilo\ampere\per\meter} \\
$A_{11}$, $A_{22}$ & \SI{6.59}{\pico\joule\per\meter} \\
$A_{12}$ & 0 \\
$A_0$ & \SI{-6.59}{\pico\joule\per\meter} \\
$a$ & \SI{1}{\nano\meter} \\
$K$ & \SI{0.15}{\mega\joule\per\meter\cubed} \\
$\aa$ & $\ee_3$ \\
$\DD$ & $D (-\ee_1 \otimes \ee_2 + \ee_2 \otimes \ee_1)$ \\
$D$ & \SI{3}{\milli\joule\per\meter\squared}
\end{tabular}
\caption{Experiment of Section~\ref{sec:numerics_skyrmion}:
Material parameters.
All values are taken from~\cite{sztfm2020}, except those of $a$ and $D$.
Here, we denote by $\{ \ee_1, \ee_2, \ee_3 \}$ the canonical basis of $\R^3$.}
\label{tab:material_parameters}
\end{table}

For the discretization, we consider a tetrahedral mesh $\T_h$ generated by Netgen
with mesh size \SI{3.36}{\nano\meter} (\num{1660} vertices and \num{4694} elements),
i.e., well below the exchange length of $\ell_{\mathrm{ex}} = \sqrt{2A_{11}/(\mu_0 \Msl{1}^2)} =$ \SI{8.61}{\nano\meter},
Here, $\mu_0 > 0$ denotes the vacuum permeability (in \si{\newton\per\ampere\squared}).

\begin{figure}[ht]
\centering
\includegraphics[width=0.4\textwidth]{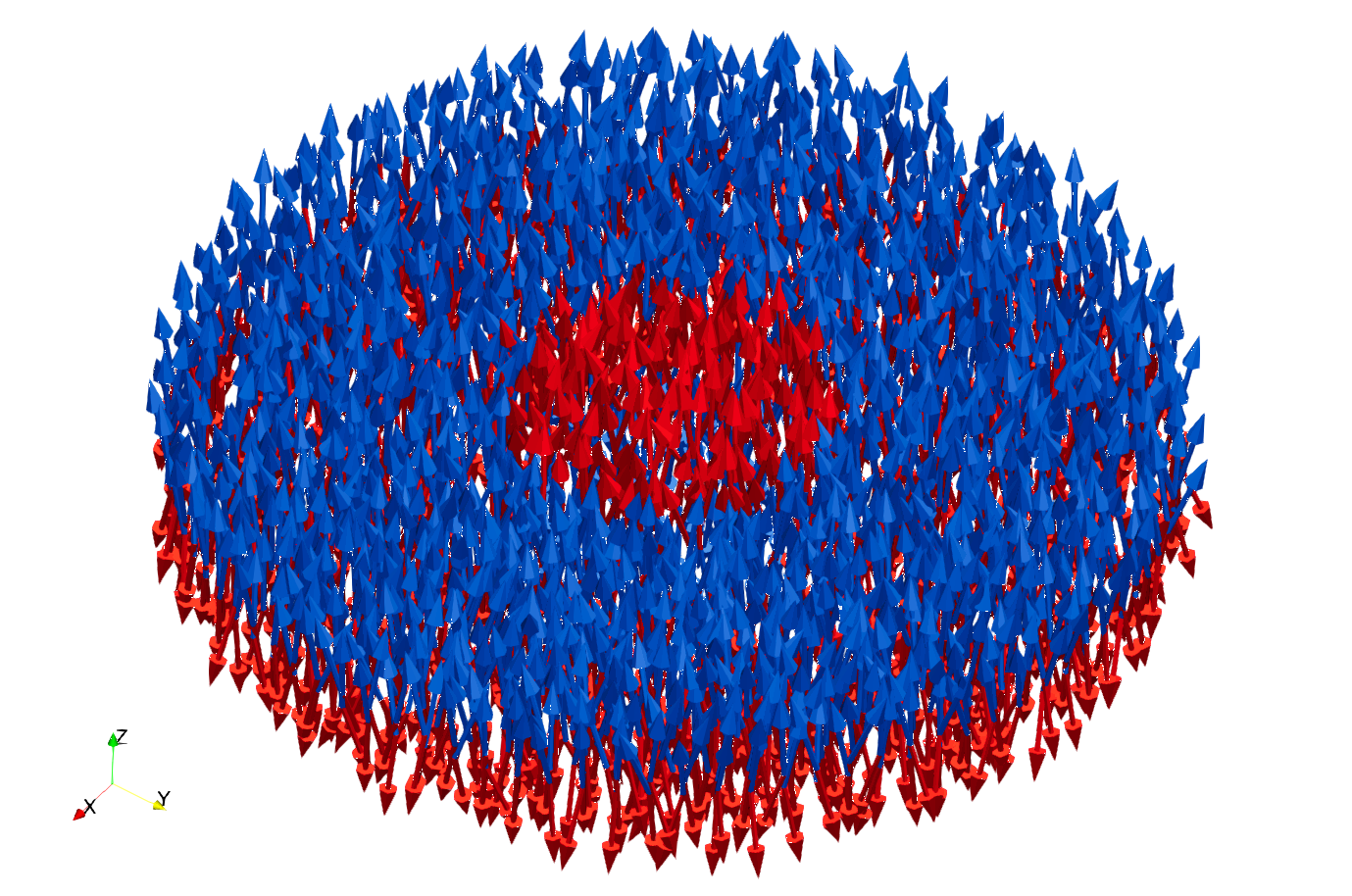}
\caption{Experiment of Section~\ref{sec:numerics_skyrmion}:
Initial guess for Algorithm~\ref{alg:decoupled}.
The initial magnetisation for $\mm_{h,1}^0$ (resp., $\mm_{h,2}^0$) is shown in red (resp., blue),
with the internal region facing in the $\ee_3$ (resp., $-\ee_3$) direction.}
\label{fig:skyrmioninitial}
\end{figure}

As an initial guess, we consider a perturbed skyrmion-like AFM state; see Figure~\ref{fig:skyrmioninitial}.
More precisely, we consider the auxiliary function
\begin{equation*}
	f_{\text{init}}(x,y,z) := \frac{1}{1 + \exp\left(20\left(\sqrt{x^2 + y^2} - 10\right)\right)}  - \frac{1}{2}
\end{equation*}
and start from the initial condition $\mm_{1,2} = (0,0,\pm f_{\text{init}})$.
This is then interpolated using the built-in Oswald-type interpolation of NGSolve
before undergoing a nodal projection, random perturbation (up to $0.3$ in each component) and another nodal projection. The value $10$ in the expression of $f_{\text{init}}$ corresponds to $\SI{10}{\nano\meter}$
and refers to the radius of the inner circle.
The decay constant $20$ makes the transition reasonably sharp before projecting.

Starting from this configuration,
we run Algorithm~\ref{alg:decoupled} (in our opinion, the best performing one in Section~\ref{sec:num_compare})
with dimensionless time-step size $\tau$ and stopping tolerance $\eps$ both equal to \num{1e-3}.

\begin{figure}[ht]
\centering
\begin{subfigure}[b]{0.4\textwidth}
\centering
\includegraphics[width=\textwidth]{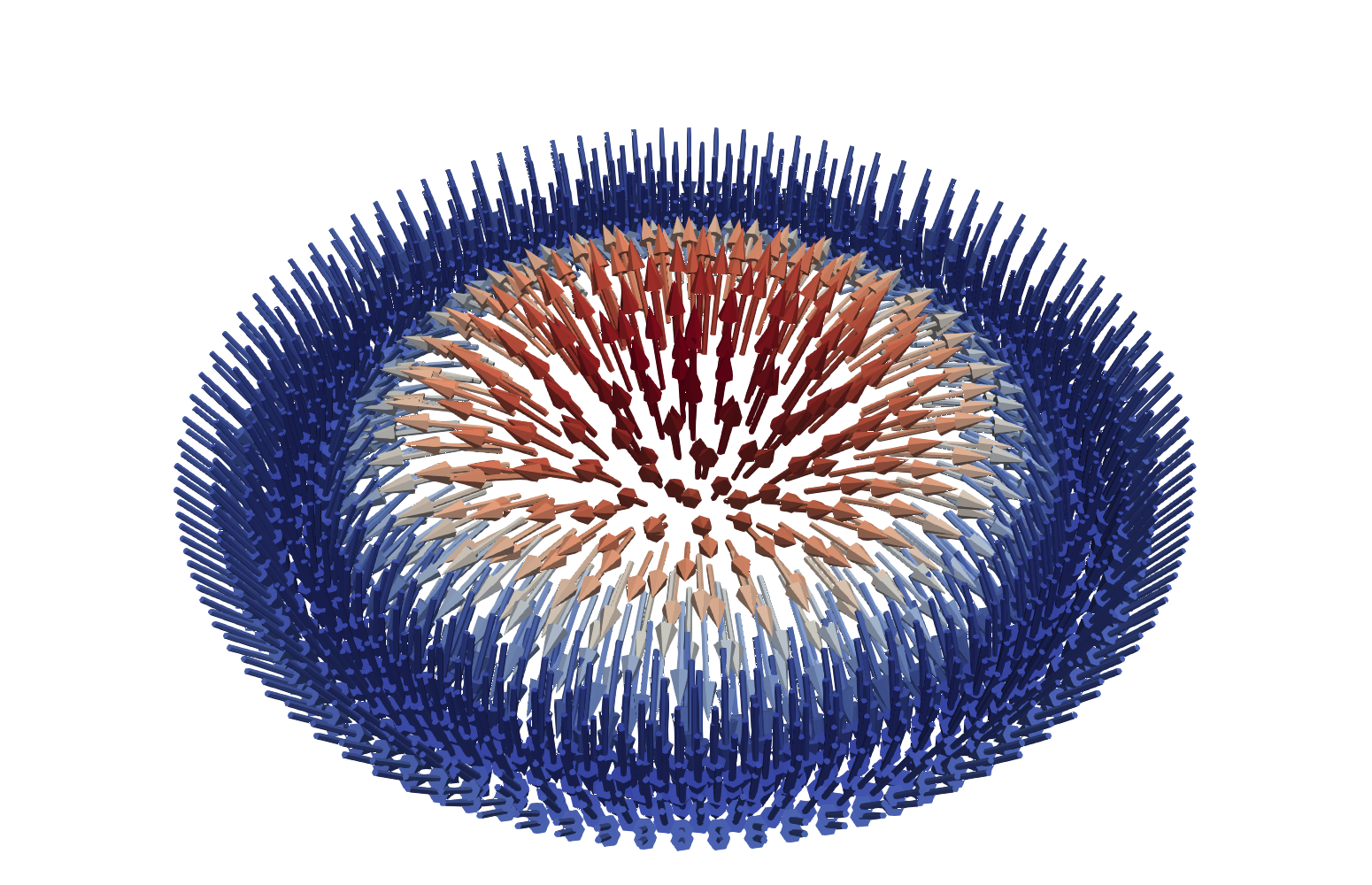}\label{fig:SkyrmionRedFinal}
\caption{$\mm_{h,1}$}
\end{subfigure}
\quad
\begin{subfigure}[b]{0.4\textwidth}
\centering
\includegraphics[width=\textwidth]{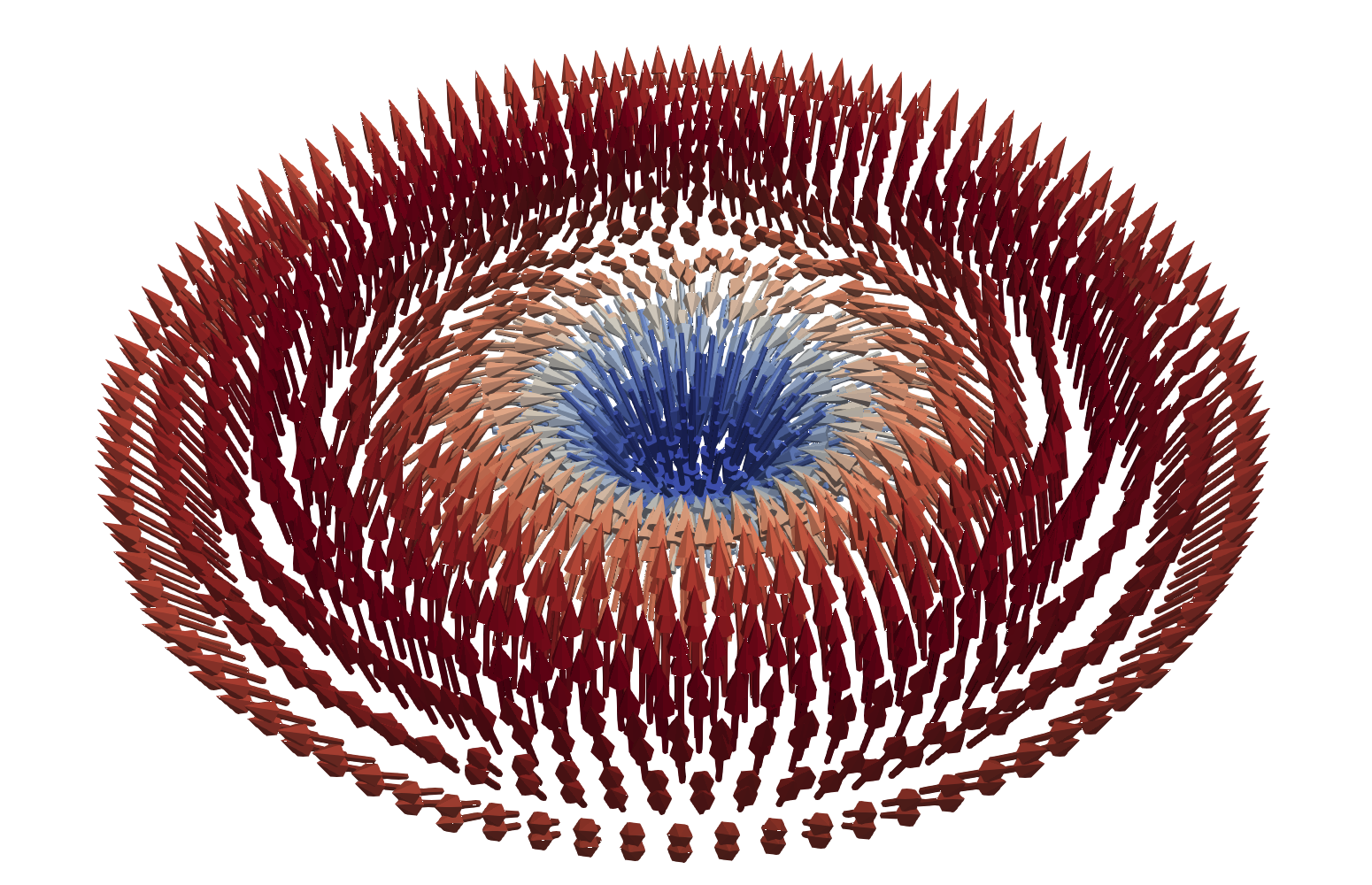}\label{fig:SkyrmionBlueFinal}
\caption{$\mm_{h,2}$}
\end{subfigure}
\caption{Experiment of Section~\ref{sec:numerics_skyrmion}:
Stable AFM configurations computed using Algorithm~\ref{alg:decoupled}
with $H^1$-metric.
In the pictures, the color scale refers to the third component of the fields,
which attains values between -1 (blue) and 1 (red).}
\label{fig:skyrmions}
\end{figure}

In Figure~\ref{fig:skyrmions},
we show the stable configurations obtained running Algorithm~\ref{alg:decoupled} with $H^1$-metric.
We see that both fields are N\'eel-type skyrmions~\cite{frc2017}
(with the cores pointing up for $\mm_{h,1}$ and down for $\mm_{h,2}$,
in line with the orientation of the field in the internal region for the corresponding initial condition),
which is typical for magnetic systems characterized by interfacial Dzyaloshinskii--Moriya interaction.
Moreover, as expected for an AFM material, we have that $\mm_{h,1} \approx -\mm_{h,2}$.

\begin{table}[ht]
\scriptsize
\begin{tabular}{c|c|c|c}
& \multicolumn{3}{c}{Algorithm~\ref{alg:decoupled}}\\
\hline
$\HHH$ & $(\LL^2(\Omega),\norm{\cdot})$ & $(\LL^2(\Omega),\norm[h]{\cdot})$ & $\HH^1(\Omega)$\\ 
\hline
\emph{energy} & \num{-4.362e5} & \num{-4.287e5} & \num{-4.256e5} \\
\emph{num.\ iter.} & \num{17552} & \num{17704} & \num{17784} \\
\emph{solve time} & 0.060 & 0.050 & 0.058 \\
$\mathrm{err}_{L^\infty}$ & \num{0.165} & \num{0.100} & \num{0.061}\\
$\mathrm{err}_{L^1}$ & \num{95.134} & \num{44.173} & \num{23.411}\\
\end{tabular}
\caption{Experiment of Section~\ref{sec:numerics_skyrmion}:
Comparison of gradient flow metrics for Algorithm~\ref{alg:decoupled}.}
\label{tab:skyrmion_minimization}
\end{table}

In Table~\ref{tab:skyrmion_minimization},
we compare the performance of the three gradient flow metrics considered in Section~\ref{sec:num_compare}.
We see that, in terms of final energy value, number of iterations, and average solve time,
the performance of the three metrics is comparable.
On the other hand, as far as the violation of the unit-length constraint is concerned,
the $H^1$-metric exhibits the best performance.

\section{Numerical approximation of the LLG system} \label{sec:llg}

In this section,
starting from Algorithm~\ref{alg:decoupled},
we introduce a fully discrete algorithm to approximate solutions of
the initial boundary value problem~\eqref{eq:llg}--\eqref{eq:llg_ibc}
for the coupled system of LLG equations
modeling the dynamics of AFM and FiM materials.
We state well-posedness, stability, and unconditional convergence
of the approximations toward a weak solution of the problem,
and present numerical experiments to show its applicability
for the simulation of the dynamics of magnetic skyrmions in AFM materials.
To make the presentation of the results concise,
all proofs are postponed to Section~\ref{sec:proofs_llg}.

\subsection{Numerical algorithm and main results} \label{sec:results_llg}

The method we propose, stated in the following algorithm,
is based on the projection-free tangent plane scheme from~\cite{ahpprs2014,bartels2016,hpprss2019}
and employs the decoupled approach of Algorithm~\ref{alg:decoupled}.
Like in the static case, the spatial discretization is based on first-order finite elements (see Section~\ref{sec:preliminaries}).

\begin{algorithm}[tangent plane scheme] \label{alg:tps}
\textbf{Discretization parameters:}
Mesh size $h>0$,
time-step size $\tau>0$.\\
\textbf{Input:}
Approximate initial condition
$(\mm_{1,h}^0, \mm_{2,h}^0) \in \S^1(\T_h)^3 \times \S^1(\T_h)^3$
such that,
for all $\ell=1,2$,
$\abs{\mm_{\ell,h}^0(z)} = 1$ for all $z \in \NN_h$.\\
\textbf{Loop:}
For all $i \in \N_0$, iterate {\rm(i)--\rm(ii)} until the stopping criterion {\rm(stop)} is met:
\begin{itemize}
\item[\rm(i)] Given $(\mm_{1,h}^i, \mm_{2,h}^i) \in \S^1(\T_h)^3 \times \S^1(\T_h)^3$,
compute $(\vv_{1,h}^i, \vv_{2,h}^i) \in \KKK_h[\mm_{1,h}^i] \times \KKK_h[\mm_{2,h}^i]$ such that,
for all $\ell=1,2$,
for all $(\pphi_{1,h} , \pphi_{2,h}) \in \KKK_h[\mm_{1,h}^i] \times \KKK_h[\mm_{2,h}^i]$,
it holds that
\begin{multline} \label{eq:tps1}
\alpha_\ell \inner[h]{\vv_{\ell,h}^i}{\pphi_{\ell,h}}
+ \inner[h]{\mm_{\ell,h}^i \times \vv_{\ell,h}^i}{\pphi_{\ell,h}}
+ \eta_\ell a_{\ell\ell} \tau \inner{\Grad\vv_{\ell,h}^i}{\Grad\pphi_{\ell,h}} \\
=
- \eta_\ell a_{\ell\ell} \inner{\Grad\mm_{\ell,h}^i}{\Grad\pphi_{\ell,h}}
- \eta_\ell a_{12} \inner{\Grad\mm_{3-\ell,h}^i}{\Grad\pphi_{\ell,h}}
+ \eta_\ell a_0 \inner{\mm_{3-\ell,h}^i}{\pphi_{\ell,h}}.
\end{multline}
\item[\rm(ii)] Define
\begin{equation} \label{eq:tps2}
\mm_{\ell,h}^{i+1} := \mm_{\ell,h}^i + \tau \vv_{\ell,h}^i
\quad \text{for all } \ell =1,2.
\end{equation}
\end{itemize}
\textbf{Output:}
Sequence of approximations $\{(\mm_{1,h}^i,\mm_{2,h}^i)\}_{i \in \N_0}$.
\end{algorithm}

Starting from approximations $\mm_{1,h}^0, \mm_{2,h}^0 \in \S^1(\T_h)^3$ of the initial conditions,
in each step of Algorithm~\ref{alg:tps},
the new approximations are computed updating the current ones using a predictor-corrector approach.
In the predictor step, \eqref{eq:tps1} are discretizations of
\begin{equation*}
\alpha_\ell \, \mmt_\ell
+ \mm_\ell \times \mmt_\ell
= \eta_\ell \, \heff{\ell}[\mm_1,\mm_2]
- \eta_\ell \, ( \heff{\ell}[\mm_1,\mm_2] \cdot \mm_\ell) \, \mm_\ell
\quad
\text{for all } \ell = 1,2,
\end{equation*}
an equivalent reformulation of~\eqref{eq:llg} that can be obtained using standard vector identities
as well as the relations $\abs{\mm_\ell} = 1$ and $\mm_\ell \cdot \mmt_\ell = 0$~\cite{aj2006}.
The discrete problems are posed in the discrete tangent space~\eqref{eq:discrete_tangent_space},
which yields a natural linearization.
Like in Algorithm~\ref{alg:decoupled},
the inhomogeneous intralattice exchange contribution is treated implicitly,
whereas the interlattice contributions are treated explicitly.
By doing this, the system of LLG equations is decoupled and one has to solve two,
independent of each other, $2N_h$-by-$2N_h$ linear systems per time-step.
The corrector step~\eqref{eq:tps2} is a simple projection-free first-order time-stepping.

In the following proposition, we show that Algorithm~\ref{alg:tps} is well-defined.

\begin{proposition}[well-posedness of Algorithm~\ref{alg:tps}] \label{prop:wellposedness}
For all $i \in \N_0$,
\eqref{eq:tps1} admits a unique solution
$(\vv_{1,h}^i, \vv_{2,h}^i) \in \KKK_h[\mm_{1,h}^i] \times \KKK_h[\mm_{2,h}^i]$.
In particular, each iteration of Algorithm~\ref{alg:tps} is well-defined.
\end{proposition}

In the following proposition, we characterize the energy behavior of Algorithm~\ref{alg:tps}.

\begin{proposition}[discrete energy law and stability of Algorithm~\ref{alg:tps}]
\label{prop:energy}
There hold the following statements:\\
{\rm(i)}
For all $i \in \N_0$,
the approximations generated by Algorithm~\ref{alg:tps} satisfy the identity
\begin{multline} \label{eq:llg_energy_law}
\E[\mm_{1,h}^{i+1},\mm_{2,h}^{i+1}]
- \E[\mm_{1,h}^i,\mm_{2,h}^i]
=
- \tau \sum_{\ell=1}^2 \frac{\alpha_\ell}{\eta_\ell} \norm[h]{\vv_{\ell,h}^i}^2
- \frac{\tau^2}{2} \sum_{\ell=1}^2 a_{\ell\ell} \norm{\Grad\vv_{\ell,h}^{i}}^2 \\
+ a_{12} \tau^2 \inner{\Grad\vv_{1,h}^{i}}{\Grad\vv_{2,h}^{i}}
- a_0 \tau^2 \inner{\vv_{1,h}^{i}}{\vv_{2,h}^{i}}.
\end{multline}
{\rm(ii)}
If $\tau < 2\max\{ \alpha_1,\alpha_2 \}/\abs{a_0}$,
for all $j \in \N$,
the approximations generated by Algorithm~\ref{alg:tps} satisfy the inequality
\begin{equation} \label{eq:llg_stability}
\sum_{\ell=1}^2 \norm[\HH^1(\Omega)]{\mm_{\ell,h}^j}^2
+ \tau \sum_{i=0}^{j-1} \sum_{\ell=1}^2 \norm[h]{\vv_{\ell,h}^i}^2
+ \tau^2 \sum_{i=0}^{j-1} \sum_{\ell=1}^2 \norm{\Grad\vv_{\ell,h}^i}^2
\le C.
\end{equation}
The constant $C>0$ depends only on the problem data
and the shape-regularity of the family of meshes.
\end{proposition}

The discrete energy law of Algorithm~\ref{alg:tps}
is an approximation of the one satisfied by weak solutions (see~\eqref{eq:weak:energy}).
The LLG-inherent energy dissipation, modulated by the damping parameters $\alpha_1$ and $\alpha_2$,
is enhanced by the dissipation coming from the second term on the right-hand side,
which is due to the implicit treatment of the homogeneous intralattice exchange contribution.
The last two terms on the right-hand side of~\eqref{eq:llg_energy_law},
in general unsigned, are perturbations arising from the explicit treatment of the
interlattice exchange contributions.

With the sequence of approximations delivered by Algorithm~\ref{alg:tps},
for $\ell=1,2$,
we define the piecewise affine time reconstruction $\mm_{\ell,h\tau} : [0,\infty) \to \S^1(\T_h)^3$ as
\begin{equation*}
\mm_{\ell,h\tau}(t) := \frac{t-t_i}{\tau}\mm_{\ell,h}^{i+1} + \frac{t_{i+1} - t}{\tau}\mm_{\ell,h}^i
\quad
\text{for all } i \in \N_0
\text{ and } t \in [t_i,t_{i+1}]
\end{equation*}
(see~\eqref{eq:reconstructions}).
In the following theorem, we state the convergence of the finite element approximations
toward a weak solution of~\eqref{eq:llg} in the sense of Definition~\ref{def:weak}.

\begin{theorem}[convergence of Algorithm~\ref{alg:tps}] \label{thm:main_llg}
Suppose that $\mm_{1,h}^0 \to \mm_1^0$ and $\mm_{2,h}^0 \to \mm_2^0$ in $\HH^1(\Omega)$ as $h \to 0$.
Then, there exist $(\mm_1, \mm_2) \in L^\infty(0,\infty;\XXX)$
and a (nonrelabeled) subsequence of $\{ (\mm_{1,h\tau} , \mm_{2,h\tau}) \}$
which converges toward $(\mm_1,\mm_2)$ as $h,\tau \to 0$.
In particular,
as $h,\tau \to 0$,
for all $\ell=1,2$
it holds that
$\mm_{\ell,h\tau} \weakstarto \mm_\ell$ in $L^{\infty}(0,\infty;\HH^1(\Omega))$.
If $a_{12} = 0$,
the limit $(\mm_{1},\mm_{2})$ is
a weak solution of~\eqref{eq:llg} in the sense of Definition~\ref{def:weak}.
\end{theorem}

Like in the stationary case,
we need to assume that $a_{12} = 0$ to be able to show that the limit toward which the finite element approximations
are converging satisfies the variational formulation~\eqref{eq:weak:variational}
(cf.\ Remark~\ref{rem:why_assumption}).
Under this assumption,
Theorem~\ref{thm:main_llg} shows existence of a weak solution to~\eqref{eq:llg}
and convergence (without rates) of the time reconstructions generated using the
snapshots computed using Algorithm~\ref{alg:tps} toward it.

\subsection{Numerical experiments} \label{sec:numerics_llg}


In this section,
we aim to show the capability of Algorithm~\ref{alg:tps} to simulate dynamic processes involving AFM materials.

\subsubsection{LLG-based energy minimization} \label{sec:num_compare_llg}

Starting from the observation that the dynamics of $\mm_1$ and $\mm_2$ governed by the system of LLG equations~\eqref{eq:llg} is dissipative,
with the energy dissipation being modulated by the damping parameters $\alpha_1$ and $\alpha_2$,
we repeat the experiment of Section~\ref{sec:num_compare},
but to compute low-energy stationary points we use Algorithm~\ref{alg:tps}
(instead of the gradient flow-based approaches of Algorithm~\ref{alg:coupled} and Algorithm~\ref{alg:decoupled}).
More precisely,
we consider the same setup and the same spatial discretization of Section~\ref{sec:num_compare}
and run Algorithm~\ref{alg:tps}
with $\eta_1=\eta_2=1$,
different damping parameters $\alpha_1 = \alpha_2 = \alpha\in\{1,1/2,1/4,1/8, 1/16 \}$,
and $\tau = 10^{-3}$,
using the constant fields $\mm_{1,h}^0 \equiv (1,0,0)$ and $\mm_{2,h}^0 \equiv (0,1,0)$ as initial condition.
The iteration is stopped when the $\alpha$-independent stopping criterion~\eqref{eq:coupled_stopping}
with $\norm[\HHH]{\cdot}^2 = \norm[h]{\cdot}^2$ and $\eps = 10^{-4}$ is satisfied.

\begin{table}[ht]
\scriptsize
\begin{tabular}{c|c|c|c|c|c|c} 
& Alg.~\ref{alg:decoupled} & \multicolumn{5}{c}{Algorithm~\ref{alg:tps}}\\
\hline
$\alpha$ & 1 & 1 & 1/2 & 1/4 & 1/8 & 1/16 \\ 
\hline
\emph{energy} & \num{-111.38} & \num{-112.26} & \num{-126.64} & \num{-160.60} & \num{-214.79} & \num{-998.44} \\
\emph{proj.\ energy err.} & \num{9.90e-11} & \num{1.03e-10} & \num{4.81e-11} & \num{8.27e-11} & \num{6.13e-11} & \num{3.64e-11} \\
\emph{num.\ iter.} & 275 & 310 & 334 & 600 & 2954 & 46698 \\
\emph{solve time} & 0.094 & 0.093 & 0.095 & 0.098 & 0.097 & 0.094 \\
$\mathrm{err}_{L^\infty}$ & 0.047 & \num{3.932e-2} & \num{9.497e-2} & \num{0.219} & \num{0.408} & \num{6.705} \\
$\mathrm{err}_{L^1}$ & \num{9.61e-02} & \num{8.019e-2} & \num{0.199} & \num{0.486} & \num{0.982} & 3.995
\end{tabular}
\caption{Experiment of Section~\ref{sec:num_compare_llg}:
Comparison of Algorithm~\ref{alg:decoupled}
(with mass-lumped $L^2$-metric and $\alpha=1$)
with Algorithm~\ref{alg:tps} (with $\alpha = 1,1/2,1/4,1/8,1/16$).}
\label{tab:constant_initial_test_llg}
\end{table}

We display the results of our computations in Table~\ref{tab:constant_initial_test_llg}.
Noting that Algorithm~\ref{alg:tps} coincides with Algorithm~\ref{alg:decoupled}
with mass-lumped $L^2$-metric
if $\eta_1 = \eta_2 = \alpha_1 = \alpha_2 = 1$
and the precession term $\inner[h]{\mm_{\ell,h}^i \times \vv_{\ell,h}^i}{\pphi_{\ell,h}}$
is omitted from~\eqref{eq:tps1},
in the first column of the table we include the results from Section~\ref{sec:num_compare}
of this instance of Algorithm~\ref{alg:decoupled}
for the sake of comparison.
We see that as $\alpha$ is lowered, the final energy is further from the expected
value of $-100$, which is due to the slower dissipation resulting in 
lengthier dynamics (larger number of iterations) and more rotations
(as the precession term is made stronger in a relative sense) before reaching the minimizing state,
thereby increasing the average length of $\mm_{h,1}$ and $\mm_{h,2}$ (as seen in the error rows).
Similarly to Table~\ref{tab:constant_initial_test} we see that after applying a nodal projection,
the energy is within 10 decimal places of $-100$, indicating that a minimizer
is still identified. As expected the average solve time is independent of $\alpha$.
For $\alpha = 1/16$ we see that the violation of the unit-length constraint
and the number of iterations are significantly larger.
We suppose that this is related to a possible instability of Algorithm~\ref{alg:tps},
since, as shown Proposition~\ref{prop:energy}(ii),
stability requires $\tau$ to be sufficiently small,
with the threshold for the time-step size being proportional to the damping parameter.
Indeed, for $\alpha=1/16$, we observe that it is sufficient to reduce $\tau$ to regain a good performance of the algorithm.

\begin{figure}[ht]
	\centering
	\begin{subfigure}[b]{0.32\textwidth}
		\centering
		\begin{tikzpicture}
			\begin{axis}[
				width = \textwidth,
				xlabel = {\scriptsize iteration},
				]
				\addplot[red, thick]	table[x=t, y=x_avg_1, col sep=comma]{DataFiles/Constant_Decoupled_L2Mass.dat};
				\addplot[blue, thick]	table[x=t, y=x_avg_2, col sep=comma]{DataFiles/Constant_Decoupled_L2Mass.dat};
			\end{axis}
		\end{tikzpicture}
		\caption{Alg.~\ref{alg:decoupled}, $\alpha=1$.}
	\end{subfigure}
	\hfill
	\begin{subfigure}[b]{0.32\textwidth}
		\centering
		\begin{tikzpicture}
			\begin{axis}[
				width = \textwidth,
				xlabel = {\scriptsize iteration},
				]
				\addplot[red, thick]	table[x=t, y=x_avg_1, col sep=comma]{DataFiles/ALPHA1.dat};
				\addplot[blue, thick]	table[x=t, y=x_avg_2, col sep=comma]{DataFiles/ALPHA1.dat};
			\end{axis}
		\end{tikzpicture}
		\caption{Alg.~\ref{alg:tps}, $\alpha=1$.}
	\end{subfigure}
	\hfill
	\begin{subfigure}[b]{0.32\textwidth}
		\centering
		\begin{tikzpicture}
			\begin{axis}[
				width = \textwidth,
				xlabel = {\scriptsize iteration},
				]
				\addplot[red, thick]	table[x=t, y=x_avg_1, col sep=comma]{DataFiles/ALPHA05.dat};
				\addplot[blue, thick]	table[x=t, y=x_avg_2, col sep=comma]{DataFiles/ALPHA05.dat};
			\end{axis}
		\end{tikzpicture}
		\caption{Alg.~\ref{alg:tps}, $\alpha=1/2$.}
	\end{subfigure}\\
	\bigskip\bigskip
	\begin{subfigure}[b]{0.32\textwidth}
		\centering
		\begin{tikzpicture}
			\begin{axis}[
				width = \textwidth,
				xlabel = {\scriptsize iteration},
				]
				\addplot[red, thick]	table[x=t, y=x_avg_1, col sep=comma]{DataFiles/ALPHA025.dat};
				\addplot[blue, thick]	table[x=t, y=x_avg_2, col sep=comma]{DataFiles/ALPHA025.dat};
			\end{axis}
		\end{tikzpicture}
		\caption{Alg.~\ref{alg:tps}, $\alpha=1/4$.}
	\end{subfigure}
	\hfill
	\begin{subfigure}[b]{0.32\textwidth}
		\centering
		\begin{tikzpicture}
			\begin{axis}[
				width = \textwidth,
				xlabel = {\scriptsize iteration},
				]
				\addplot[red, thick]	table[x=t, y=x_avg_1, col sep=comma]{DataFiles/ALPHA0125.dat};
				\addplot[blue, thick]	table[x=t, y=x_avg_2, col sep=comma]{DataFiles/ALPHA0125.dat};
			\end{axis}
		\end{tikzpicture}
		\caption{Alg.~\ref{alg:tps}, $\alpha=1/8$.}
	\end{subfigure}
	\hfill
	\begin{subfigure}[b]{0.32\textwidth}
		\centering
		\begin{tikzpicture}
			\begin{axis}[
				width = \textwidth,
				xlabel = {\scriptsize iteration},
				]
				\addplot[red, thick]	table[x=t, y=x_avg_1, col sep=comma]{DataFiles/ALPHA00625.dat};
				\addplot[blue, thick]	table[x=t, y=x_avg_2, col sep=comma]{DataFiles/ALPHA00625.dat};
			\end{axis}
		\end{tikzpicture}
		\caption{Alg.~\ref{alg:tps}, $\alpha=1/16$.}
	\end{subfigure}
	\caption{Experiment of Section~\ref{sec:numerics_llg}:
		Evolution of $\langle \mm_1(t) \cdot\ee_1 \rangle$ (red)
		and $\langle \mm_2(t) \cdot\ee_1 \rangle$ (blue).
		(a) Algorithm~\ref{alg:decoupled} with mass-lumped $L^2$-metric and $\alpha = 1$.
		(b)--(f) Algorithm~\ref{alg:tps} with $\alpha = 1, 1/2, 1/4, 1/8, 1/16$.
	}
	\label{fig:constant_llg}
\end{figure}
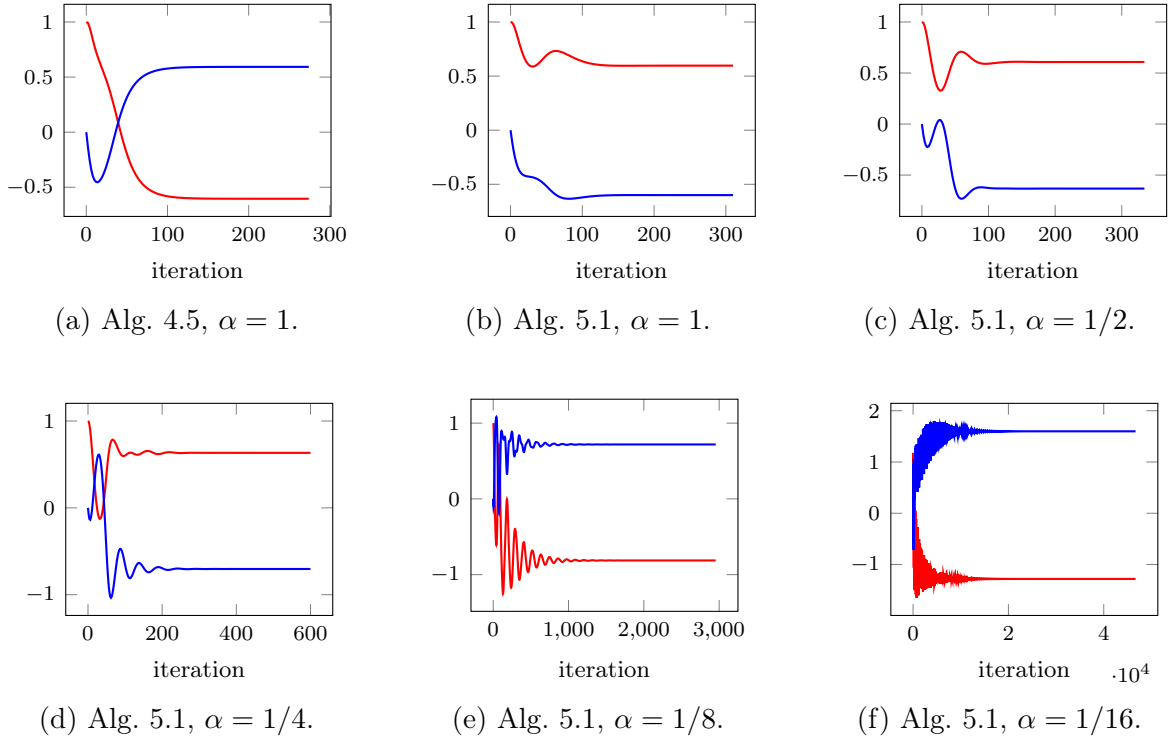

The fact that lowering the value of $\alpha$ results in less dissipation
can also be seen in Figure~\ref{fig:constant_llg},
where we display
the evolution of the average first component of both fields,
i.e., $\langle \mm_{\ell}(t) \cdot\ee_1 \rangle = \abs{\Omega}^{-1}\int_{\Omega} \mm_{\ell}(t) \cdot\ee_1$
for all $\ell=1,2$, for all cases.
Interestingly,
we see that the gradient flow dynamics and the LLG dynamics for $\alpha=1/8,1/16$
return as approximate stationary point an approximation of the minimizer $(\mm_1^-,\mm_2^-) \equiv (-\aa, \aa)$,
whereas the LLG dynamics for $\alpha=1,1/2,1/4$ return an approximation of $(\mm_1^+,\mm_2^+) \equiv (\aa, -\aa)$.
This is not surprising, since different dynamics can result in convergence to different stationary points,
even with the same initial condition.
As far as the LLG dynamics is concerned, we see that the oscillations of the average first components
increase as $\alpha$ is lowered, which can be explained by the greater relative weight of the precessional term
on the right-hand side of the LLG equation for smaller values of $\alpha$.

\subsubsection{Skyrmion dynamics} \label{sec:numerics_skyrmion_llg}

Inspired the experiment in~\cite[Section~4.3]{hpprss2019}, we simulate the dynamics
of isolated magnetic skyrmions in an AFM nanodisk in response to an applied field pulse.

The setup (domain, energy, and material parameters) is the same as in Section~\ref{sec:numerics_skyrmion},
which we complete with the additional parameters needed for the dynamic case,
i.e., the rescaled gyromagnetic ratios
$\gamma_1 = \gamma_2 = \gamma_0 \approx$ \SI{2.21e5}{\meter\per\ampere\per\second}
and the Gilbert damping parameters $\alpha_1 = \alpha_2 =$ \num{5e-3}
(see~\eqref{eq:LLG} below).
Given the same spatial discretization (mesh) as in Section~\ref{sec:numerics_skyrmion},
as initial conditions $\mm_{1,h}^0$ and $\mm_{2,h}^0$ for Algorithm~\ref{alg:tps},
we consider the nodal projections of the N\'eel-type skyrmions shown in Figure~\ref{fig:skyrmions}.
Moreover, for the time discretization,
we use of a constant time-step size of \SI{2}{\femto\second}.

\begin{figure}[ht]
\centering
\begin{subfigure}[b]{0.47\textwidth}
\centering
\begin{tikzpicture}
\begin{axis}[
width = \textwidth,
xlabel={\scriptsize $t$ (\si{\nano\second})},
xmin=0,
xmax=1,
ymin=-10,
ymax=70,
xtick={0,0.150,1},
xticklabels={0,0.150,1},
ytick={0,60},
yticklabels={0,$H_{\mathrm{max}}$},
]
\addplot[blue,thick] coordinates {(0,0) (0.040,60)};
\addlegendentry{\tiny $H(t)$}
\addplot[blue,thick] coordinates {(0.040,60) (0.110,60)};
\addplot[blue,thick] coordinates {(0.110,60) (0.150,0)};
\addplot[blue,thick] coordinates {(0.150,0) (1,0)};
\addplot[dashed] coordinates {(0.040,-10) (0.040,60)};
\addplot[dashed] coordinates {(0.110,-10) (0.110,60)};
\end{axis}
\end{tikzpicture}
\caption{Applied pulse field.}
\end{subfigure}
\hfill
\begin{subfigure}[b]{0.47\textwidth}
\centering
\begin{tikzpicture}
	\begin{axis}[
		width = \textwidth,
		xlabel = {\scriptsize $t$ (\si{\nano\second})},
		xmin = 0,
		xmax = 1,
		]
		\addplot[teal, thick]	table[x=t, y=x_avg_total, col sep=comma]{DataFiles/SkyrmionLLG.dat};
	\end{axis}
\end{tikzpicture}
\caption{$\langle \mm(t) \cdot\ee_1 \rangle$.}
\end{subfigure}
\caption{Experiment of Section~\ref{sec:numerics_skyrmion_llg}:
(a) Structure of the applied pulse field.
(b) Time evolution of $\langle \mm(t) \cdot\ee_1 \rangle$.
}
\label{fig:skyrmion_llg}
\end{figure}
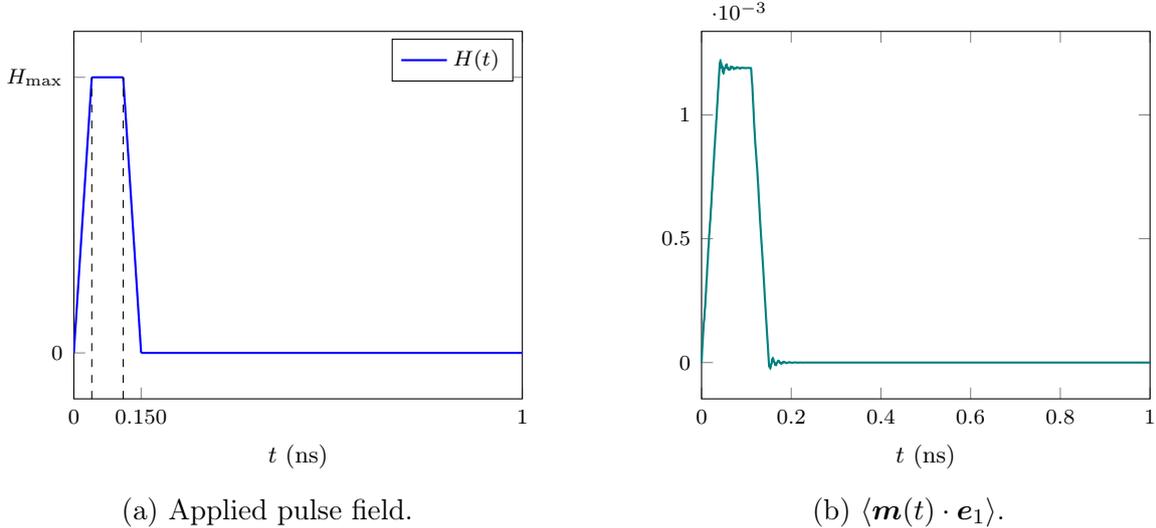

Starting from this configuration, we perturb the system from its equilibrium by applying
an in-plane pulse field of the form $\Hext(t) = (H(t),0,0)$ of maximum intensity
$\mu_0 H_{\mathrm{max}} =$ \SI{100}{\milli\tesla} for \SI{150}{ps}; see Figure~\ref{fig:skyrmion_llg}(a).
Then, we turn off the applied external field, i.e., $\Hext(t) \equiv (0,0,0)$, and let the system relax to equilibrium.
The overall simulation time is \SI{1}{\nano\second}.

In Figure~\ref{fig:skyrmion_llg}(b),
we show the time evolution of the average first component of the total magnetization $\mm=\mm_1+\mm_2$.
We see a perfect match between the applied pulse field and the total magnetization.
When the field is turned off, the state immediately comes back to the initial configuration, which confirms its stability.

\section{Proofs} \label{sec:proofs}

In this section, we collect the proof of all results presented in the paper.

\subsection{Static problem} \label{sec:proofs_min}

We start with showing the weak sequential lower semicontinuity of the energy functional.

\begin{proposition} \label{prop:wlsc_energy}
The energy functional~\eqref{eq:energy} is weakly sequentially lower semicontinuous in
$\HH^1(\Omega) \times \HH^1(\Omega)$, i.e.,
if $\{ (\mm_{1,k},\mm_{2,k}) \}_{k \in \N} \subset \HH^1(\Omega)\times \HH^1(\Omega)$
and $(\mm_1,\mm_2) \in \HH^1(\Omega)\times \HH^1(\Omega)$
are such that
$(\mm_{1,k},\mm_{2,k}) \weakto (\mm_1,\mm_2)$ in $\HH^1(\Omega)\times \HH^1(\Omega)$ as $k \to \infty$,
then $\E[\mm_1,\mm_2] \le \liminf_{k \to \infty} \E[\mm_{1,k},\mm_{2,k}]$.
\end{proposition}

The result is a special case of the following lemma.

\begin{lemma} \label{lem:wlsc_energy}
Let $V$ and $H$ two Hilbert spaces such that $V \subset H$ with compact inclusion.
Let $a: V \times V \to \R$ be a continuous bilinear form satisfying a so-called G\r{a}rding inequality,
i.e.,
there exists $C_1 >0$ and $C_2 \in \R$ such that
\begin{equation} \label{eq:garding}
a(v,v) \ge C_1 \norm[V]{v}^2 - C_2 \norm[H]{v}^2
\quad
\text{for all } v \in V.
\end{equation}
Then, the quadratic functional $\mathcal{J} : V \to \R$ defined by $\mathcal{J}[v]:= a(v,v)$ for all $v \in V$
is weakly sequentially lower semicontinuous in $V$,
i.e.,
if $\{ v_k \}_{k \in \N} \subset V$ and $v \in V$ are such that
$v_k \weakto v$ in $V$ as $k \to \infty$,
then $\mathcal{J}[v] \le \liminf_{k \to \infty} \mathcal{J}[v_k]$.
\end{lemma}

\begin{proof}
Let $\{ v_k \}_{k \in \N} \subset V$ and $v \in V$ be such that $v_k \weakto v$ in $V$ as $k \to \infty$.
From the compact inclusion $V \subset H$, it follows that $v_k \to v$ in $H$.
Using~\eqref{eq:garding}, we see that
\begin{equation*}
C_1 \norm[V]{v - v_k}^2 - C_2 \norm[H]{v - v_k}^2 \le a(v-v_k,v-v_k)
= a(v,v) - a(v_k,v) - a(v,v_k) + a(v_k,v_k).
\end{equation*}
We now take the liminf as $k \to \infty$ of this inequality.
For the left-hand side we have that
\begin{equation*}
\liminf_{k \to \infty} \left( C_1 \norm[V]{v - v_k}^2 - C_2 \norm[H]{v - v_k}^2 \right) \ge 0.
\end{equation*}
For the right-hand side, noting that $v_k \weakto v$ in $V$ implies that $a(v_k,v) \to a(v,v)$ as $k \to \infty$,
we have that
\begin{equation*}
\begin{split}
\liminf_{k \to \infty} \left[ a(v,v) - a(v_k,v) - a(v,v_k) + a(v_k,v_k) \right]
& = - a(v,v) + \liminf_{k \to \infty} a(v_k,v_k) \\
& = - \mathcal{J}[v] + \liminf_{k \to \infty} \mathcal{J}[v_k].
\end{split}
\end{equation*}
This shows that $\mathcal{J}[v] \le \liminf_{k \to \infty} \mathcal{J}[v_k]$
and thus concludes the proof.
\end{proof}

We now prove Theorem~\ref{thm:gamma_convergence} establishing the $\Gamma$-convergence
of our finite element discretization.

\begin{proof}[Proof of Theorem~\ref{thm:gamma_convergence}]
Part~(i) of the theorem immediately follows from the 
weak sequential lower semicontinuity of the energy functional established in
Proposition~\ref{prop:wlsc_energy}.

To show part~(ii), let $(\mm_1,\mm_2) \in \XXX$ be arbitrary.
Since $\CC^\infty(\overline{\Omega};\sphere)$ is dense in $\HH^1(\Omega;\sphere)$
(see \cite[Theorem~III.6.2]{struwe2008}),
for all $k \in \N$
there exists $(\mm_{1,k},\mm_{2,k}) \in \CC^\infty(\overline{\Omega};\sphere) \times \CC^\infty(\overline{\Omega};\sphere)$ such that
$\norm[\HH^1(\Omega)]{\mm_\ell - \mm_{\ell,k}} \le 1/k$ for all $\ell=1,2$.

Let $\eps>0$.
The above convergence guarantees the existence of $k \in \N$
such that $\norm[\HH^1(\Omega)]{\mm_\ell - \mm_{\ell,k}} \le \eps/2$.
Define $\mm_{\ell,k,h} := \interp[\mm_{\ell,k}]$ for all $\ell=1,2$.
By construction, for all $\ell=1,2$,
$\abs{\mm_{\ell,k,h}(z)} = 1$ for all $z \in \NN_h$
and $0 = \norm[L^1(\Omega)]{\interp[\abs{\mm_{\ell,k,h}}^2]-1} \le \delta$
for all $\delta>0$.
Hence, $(\mm_{1,k,h},\mm_{2,k,h})$ belongs to $\XXX_{h,\delta}$ for all $\delta>0$.
Moreover,
a classical interpolation estimate yields that
$\norm[\HH^1(\Omega)]{\mm_{\ell,k} - \mm_{\ell,k,h}}
\le C h \norm{D^2\mm_{\ell,k}}$.
Therefore, we have that
$\norm[\HH^1(\Omega)]{\mm_{\ell,k} - \mm_{\ell,k,h}}
\le \eps/2$ if $h$ is chosen sufficiently small.
Using the triangle inequality, we thus obtain that
$\norm[\HH^1(\Omega)]{\mm_\ell - \mm_{\ell,k,h}} \le \eps$.
Since $\eps>0$ was arbitrary,
this shows that the sequence
$\{ (\mm_{1,h,\delta} , \mm_{2,h,\delta}) \}$
defined by
$(\mm_{1,h,\delta} , \mm_{2,h,\delta}) := ((\mm_{1,k,h},\mm_{2,k,h})) \in \XXX_{h,\delta}$
satisfies the desired convergence property toward $(\mm_1,\mm_2)$
as $h,\delta \to 0$
(note that our construction is independent of $\delta$, so the limit $\delta \to 0$ is trivial).
This implies also that $\E_{h,\delta}[\mm_{1,h,\delta} , \mm_{2,h,\delta}] \to \E[\mm_1 , \mm_2]$
as $h,\delta \to 0$
and concludes the proof.
\end{proof}

In view of the analysis of the discrete gradient flows presented in Section~\ref{sec:static},
we now introduce the following algorithm.

\begin{algorithm}[general discrete gradient flow] \label{alg:general}
\textbf{Discretization parameters:}
Mesh size $h>0$,
time-step size $\tau>0$, tolerance $\eps>0$,
parameters $0 \le \theta_1, \theta_2, \theta_3 \le 1$.\\
\textbf{Input:}
Initial guess
$(\mm_{1,h}^0, \mm_{2,h}^0) \in \S^1(\T_h)^3 \times \S^1(\T_h)^3$ such that,
for all $\ell=1,2$,
$\abs{\mm_{\ell,h}^0(z)} = 1$ for all $z \in \NN_h$.\\
\textbf{Loop:}
For all $i \in \N_0$, iterate {\rm(i)--\rm(ii)} until the stopping criterion {\rm(stop)} is met:
\begin{itemize}
\item[\rm(i)] Given $(\mm_{1,h}^i, \mm_{2,h}^i) \in \S^1(\T_h)^3 \times \S^1(\T_h)^3$,
compute $(\vv_{1,h}^i, \vv_{2,h}^i) \in \KKK_h[\mm_{1,h}^i] \times \KKK_h[\mm_{2,h}^i]$ such that,
for all $(\pphi_{1,h} , \pphi_{2,h}) \in \KKK_h[\mm_{1,h}^i] \times \KKK_h[\mm_{2,h}^i]$
and $\ell=1,2$, it holds that
\begin{multline} \label{eq:general1}
\inner[\HHH]{\vv_{\ell,h}^i}{\pphi_{\ell,h}}
+ a_{\ell\ell} \theta_1 \tau \inner{\Grad\vv_{\ell,h}^i}{\Grad\pphi_{\ell,h}}
+ a_{12}  \theta_2 \tau \inner{\Grad\vv_{3-\ell,h}^i}{\Grad\pphi_{\ell,h}}
- a_0  \theta_3 \tau \inner{\vv_{3-\ell,h}^i}{\pphi_{\ell,h}} \\
=
- a_{\ell\ell} \inner{\Grad\mm_{\ell,h}^i}{\Grad\pphi_{\ell,h}}
- a_{12} \inner{\Grad\mm_{3-\ell,h}^i}{\Grad\pphi_{\ell,h}}
+ a_0 \inner{\mm_{3-\ell,h}^i}{\pphi_{\ell,h}}.
\end{multline}
\item[\rm(ii)] Define
\begin{equation} \label{eq:general2}
\mm_{\ell,h}^{i+1} := \mm_{\ell,h}^i + \tau \vv_{\ell,h}^i
\quad \text{for all } \ell =1,2.
\end{equation}
\item[\rm(stop)]
Stop iterating {\rm(i)--(ii)} if $(\vv_{1,h}^i, \vv_{2,h}^i) \in \KKK_h[\mm_{1,h}^i] \times \KKK_h[\mm_{2,h}^i]$ satisfies
\begin{equation} \label{eq:general_stopping}
\sum_{\ell=1}^2 \left( \norm[\HHH]{\vv_{\ell,h}^i}^2 + \tau \norm{\Grad\vv_{\ell,h}^i}^2 \right)
\le \eps^2 \abs{\Omega}.
\end{equation}
\end{itemize}
\textbf{Output:}
If $i^* \in \N_0$ denotes the smallest integer satisfying the stopping criterion~\eqref{eq:general_stopping},
define the approximate stationary point
$(\mm_{1,h},\mm_{2,h}) := (\mm_{1,h}^{i^*},\mm_{2,h}^{i^*})$.
\end{algorithm}

The parameters $0 \le \theta_1, \theta_2, \theta_3 \le 1$ modulates the `degree of implicitness'
in the treatment of the three contributions of the energy.
It is easy to see that
Algorithm~\ref{alg:coupled} and Algorithm~\ref{alg:decoupled} are special instances
of Algorithm~\ref{alg:general},
where $\theta_1 = 1$ (backward Euler) and $ \theta_2 =  \theta_3 = 1/2$ (Crank--Nicolson) in Algorithm~\ref{alg:coupled},
whereas $\theta_1 = 1$ (backward Euler) and $ \theta_2 =  \theta_3 = 0$ (forward Euler) in Algorithm~\ref{alg:decoupled}.

For ease of presentation, in Section~\ref{sec:static},
we have decided not to present Algorithm~\ref{alg:general} in its full generality,
but we have restricted ourselves to two of its instances.
This has been motivated by the following two reasons:
First, we believe that the two proposed cases are the most relevant in practical computations.
Second, the properties and the analysis of the algorithm for general $\theta_1, \theta_2, \theta_3$
resemble the ones of the two presented prototypical cases
(excluding the combinations involving values $\theta_1, \theta_2 < 1/2$,
which require severe restrictions of the form $\tau = \mathcal{O}(h^2)$ for stability
and therefore have been ignored).

In the following proposition,
we show well-posedness of each iteration of Algorithm~\ref{alg:general}.

\begin{proposition} \label{prop:wellposedness_general}
Suppose that $\theta_1$, $\theta_2$, $\theta_3$, and $\tau$ satisfy the following conditions:
\begin{equation} \label{eq:inequalities_wellposedness}
\theta_1 > 0,
\quad
a_{11}a_{22}\theta_1^2 > a_{12}^2 \theta_2^2,
\quad
\text{and}
\quad
\cH^2 \abs{a_0} \theta_3 \tau < 1,
\end{equation}
where $a_{11}$, $a_{22}$, $a_{12}$, and $a_0$ are the coefficients in~\eqref{eq:energy},
whereas $\cH$ is the constant in~\eqref{eq:metric}.
Then, for all $i \in \N_0$, \eqref{eq:general1} admits a unique solution
$(\vv_{1,h}^i, \vv_{2,h}^i) \in \KKK_h[\mm_{1,h}^i] \times \KKK_h[\mm_{2,h}^i]$.
\end{proposition}

\begin{proof}
Let $i \in \N_0$.
The sum of the left-hand sides of~\eqref{eq:general1} for $\ell=1,2$ yields a bilinear form
$b_i : (\KKK_h[\mm_{1,h}^i] \times \KKK_h[\mm_{2,h}^i]) \times (\KKK_h[\mm_{1,h}^i] \times \KKK_h[\mm_{2,h}^i]) \to \R$,
which is defined by
\begin{equation*}
\begin{split}
b_i ( (\ppsi_{1,h},\ppsi_{2,h}) , (\pphi_{1,h},\pphi_{2,h}) )
& =
\inner[\HHH]{\ppsi_{1,h}}{\pphi_{1,h}}
+ \inner[\HHH]{\ppsi_{2,h}}{\pphi_{2,h}} \\
& \quad
+ a_{11} \theta_1 \tau \inner{\Grad\ppsi_{1,h}}{\Grad\pphi_{1,h}}
+ a_{22} \theta_1 \tau \inner{\Grad\ppsi_{2,h}}{\Grad\pphi_{2,h}} \\
& \quad
+ a_{12}  \theta_2 \tau \inner{\Grad\ppsi_{2,h}}{\Grad\pphi_{1,h}}
+ a_{12}  \theta_2 \tau \inner{\Grad\ppsi_{1,h}}{\Grad\pphi_{2,h}} \\
& \quad
- a_0  \theta_3 \tau \inner{\ppsi_{2,h}}{\pphi_{1,h}}
- a_0  \theta_3 \tau \inner{\ppsi_{1,h}}{\pphi_{2,h}}
\end{split}
\end{equation*}
for all $(\ppsi_{1,h},\ppsi_{2,h}), (\pphi_{1,h},\pphi_{2,h}) \in \KKK_h[\mm_{1,h}^i] \times \KKK_h[\mm_{2,h}^i]$.
Owing to the second inequality in~\eqref{eq:metric}, the bilinear form is bounded with respect to the $H^1$-norm.
To show coercivity, for an arbitrary $(\pphi_{1,h},\pphi_{2,h}) \in \KKK_h[\mm_{1,h}^i] \times \KKK_h[\mm_{2,h}^i]$, we first compute
\begin{equation*}
\begin{split}
b_i ( (\pphi_{1,h},\pphi_{2,h}) , (\pphi_{1,h},\pphi_{2,h}) )
& =
\norm[\HHH]{\pphi_{1,h}}^2
+ \norm[\HHH]{\pphi_{2,h}}^2
- 2 a_0  \theta_3 \tau \inner{\pphi_{2,h}}{\pphi_{1,h}} \\
& \quad
+ a_{11} \theta_1 \tau \norm{\Grad\pphi_{1,h}}^2
+ a_{22} \theta_1 \tau \norm{\Grad\pphi_{2,h}}^2 \\
& \quad
+ 2 a_{12}  \theta_2 \tau \inner{\Grad\pphi_{2,h}}{\Grad\pphi_{1,h}}.
\end{split}
\end{equation*}
The terms involving the gradients of $(\pphi_{1,h},\pphi_{2,h})$ make up a quadratic form,
which is positive definite if and only if the underlying 2-by-2 matrix is positive definite,
which is true if and only if the first two inequalities in~\eqref{eq:inequalities_wellposedness} hold.

Thanks to~\eqref{eq:metric}, it holds that
\begin{equation*}
\norm[\HHH]{\pphi_{1,h}}^2
+ \norm[\HHH]{\pphi_{2,h}}^2
- 2 a_0  \theta_3 \tau \inner{\pphi_{2,h}}{\pphi_{1,h}}
\ge (\cH^{-2} - \abs{a_0} \theta_3 \tau) \left(\norm{\pphi_{1,h}}^2 + \norm{\pphi_{2,h}}^2 \right).
\end{equation*}
This shows that the $L^2$-part of the bilinear form is coercive
if the third inequality in~\eqref{eq:inequalities_wellposedness} holds.

Hence, we conclude that the bilinear form $b_i(\cdot,\cdot)$ is coercive with respect to the $H^1$-norm
if~\eqref{eq:inequalities_wellposedness} is satisfied.
Observing that the sum over $\ell=1,2$ of the right-hand sides of~\eqref{eq:general1}
defines a bounded linear form,
well-posedness of~\eqref{eq:general1} then follows from the Lax--Milgram theorem.
\end{proof}

In the following proposition,
we establish the discrete energy law satisfied by the approximations generated by Algorithm~\ref{alg:general}.

\begin{proposition} \label{prop:energy_law_general}
Let $\theta_1$, $\theta_2$, $\theta_3$, and $\tau$ satisfy the assumptions of Proposition~\ref{prop:wellposedness_general}.
For all $i \in \N_0$, the iterates of Algorithm~\ref{alg:general} satisfy
\begin{multline} \label{eq:general_energy_law}
\E[\mm_{1,h}^{i+1},\mm_{2,h}^{i+1}]
- \E[\mm_{1,h}^i,\mm_{2,h}^i]
=
- \tau \sum_{\ell=1}^2 \norm[\HHH]{\vv_{\ell,h}^i}^2
- \frac{(2\theta_1-1)}{2} \tau^2 \sum_{\ell=1}^2 a_{\ell\ell} \norm{\Grad\vv_{\ell,h}^i}^2 \\
- a_{12} (2\theta_2-1) \tau^2 \inner{\Grad\vv_{1,h}^i}{\Grad\vv_{2,h}^i}
+ a_0 (2\theta_3-1) \tau^2 \inner{\vv_{1,h}^i}{\vv_{2,h}^i}.
\end{multline}
Suppose that $\theta_1$, $\theta_2$, $\theta_3$, and $\tau$ satisfy also the following conditions:
\begin{equation} \label{eq:inequalities_decay}
\theta_1 \ge 1/2,
\quad
a_{11}a_{22}(2\theta_1-1)^2 \ge a_{12}^2 (2\theta_2-1)^2,
\quad
\text{and }
\cH^2 \abs{a_0} \abs{2 \theta_3 - 1} \tau \le 2.
\end{equation}
Then, the sequence of energies generated by Algorithm~\ref{alg:general} is monotonically decreasing,
i.e., it holds that $\E[\mm_{1,h}^{i+1},\mm_{2,h}^{i+1}] \le \E[\mm_{1,h}^i,\mm_{2,h}^i]$
for all $i \in \N_0$.
\end{proposition}

\begin{proof}
Let $i \in \N_0$.
Testing~\eqref{eq:general1} with $\pphi_{\ell,h} = \vv_{\ell,h}^i \in \KKK_h[\mm_{\ell,h}^i]$
for $\ell=1,2$ and summing the resulting equations,
we obtain that
\begin{equation} \label{eq_aux_stability}
\begin{split}
& \sum_{\ell=1}^2 \left(\norm[\HHH]{\vv_{\ell,h}^i}^2
+ a_{\ell\ell} \theta_1 \tau \norm{\Grad\vv_{\ell,h}^i}^2 \right)
+ 2 a_{12} \theta_2 \tau \inner{\Grad\vv_{1,h}^i}{\Grad\vv_{2,h}^i}
- 2 a_0 \theta_3 \tau \inner{\vv_{1,h}^i}{\vv_{2,h}^i} \\
& \quad = 
\sum_{\ell=1}^2 \left( - a_{\ell\ell} \inner{\Grad\mm_{\ell,h}^i}{\Grad\vv_{\ell,h}^i}
- a_{12} \inner{\Grad\mm_{3-\ell,h}^i}{\Grad\vv_{\ell,h}^i}
+ a_0 \inner{\mm_{3-\ell,h}^i}{\vv_{\ell,h}^i} \right).
\end{split}
\end{equation}
It follows that
\begin{equation*}
\begin{split}
& \E[\mm_{1,h}^{i+1},\mm_{2,h}^{i+1}] \\
& \quad \stackrel{\eqref{eq:general2}}{=}
\E[\mm_{1,h}^i,\mm_{2,h}^i]
+ \frac{1}{2} \tau \sum_{\ell=1}^2 a_{\ell\ell} \big( 2 \inner{\Grad\mm_{\ell,h}^i}{\Grad\vv_{\ell,h}^i}
+ \tau \norm{\Grad\vv_{\ell,h}^i}^2  \big) \\
& \qquad\quad
+ a_{12} \tau \Big( \inner{\Grad\mm_{1,h}^i}{\Grad\vv_{2,h}^i}
+ \inner{\Grad\mm_{2,h}^i}{\Grad\vv_{1,h}^i}
+ \tau \inner{\Grad\vv_{1,h}^i}{\Grad\vv_{2,h}^i} \big) \\
& \qquad\quad
- a_0 \tau \big( \inner{\mm_{1,h}^i}{\vv_{2,h}^i}
+ \inner{\mm_{2,h}^i}{\vv_{1,h}^i}
+ \tau \inner{\vv_{1,h}^i}{\vv_{2,h}^i} \Big) \\
& \quad \stackrel{\eqref{eq_aux_stability}}{=}
\E[\mm_{1,h}^i,\mm_{2,h}^i]
- \tau \sum_{\ell=1}^2 \norm[\HHH]{\vv_{\ell,h}^i}^2
- \frac{(2\theta_1-1)}{2} \tau^2 \sum_{\ell=1}^2 a_{\ell\ell} \norm{\Grad\vv_{\ell,h}^i}^2 \\
& \qquad\quad
- a_{12} (2\theta_2-1) \tau^2 \inner{\Grad\vv_{1,h}^i}{\Grad\vv_{2,h}^i}
+ a_0 (2\theta_3-1) \tau^2 \inner{\vv_{1,h}^i}{\vv_{2,h}^i},
\end{split}
\end{equation*}
which is~\eqref{eq:general_energy_law}.
Arguing as in the proof of Proposition~\ref{prop:wellposedness_general},
it is easy to see the right-hand side of~\eqref{eq:general_energy_law} is
nonpositive if the inequalities in~\eqref{eq:inequalities_decay} are satisfied.
This shows that the sequence of energies generated by the algorithm is monotonically decreasing and concludes the proof.
\end{proof}

In the following lemma,
we prove two auxiliary estimates,
which will be useful in the proof of convergence of Algorithm~\ref{alg:general}.

\begin{lemma} \label{lem:projection_free_lemma}
For all $\ell =1,2$, for all $j \in \N$,
the iterates of Algorithm~\ref{alg:general} satisfy
\begin{align}
\label{eq:general_constraint}
\cT^{-1} \norm[L^1(\Omega)]{\interp[\abs{\mm_{\ell,h}^j}^2]-1}
& \le \tau^2 \sum_{i=0}^{j-1} \norm{\vv_{\ell,h}^i}^2 \\
\label{eq:general_L2_boundedness}
\cT^{-1} \norm{\mm_{\ell,h}^j}^2
& \le \abs{\Omega} + \tau^2 \sum_{i=0}^{j-1} \norm{\vv_{\ell,h}^i}^2,
\end{align}
where $\cT>0$ depends only on the shape-regularity of the family of meshes.
\end{lemma}

\begin{proof}
We follow~\cite{bartels2016}.
Let $\ell =1,2$ and $j \in \N$.
For all $i=0,\dots,j-1$,
from~\eqref{eq:general2},
since $\vv_{\ell,h}^i \in \KKK_h[\mm_{\ell,h}^i]$,
we deduce that
$\abs{\mm_{\ell,h}^{i+1}(z)}^2 = \abs{\mm_{\ell,h}^i(z)}^2 + \tau^2\abs{\vv_{\ell,h}^i(z)}^2$
for all $z \in \NN_h$.
Iterating in $i$ and using that $\abs{\mm_{\ell,h}^0(z)}=1$ for all $z \in \NN_h$, we obtain that
\begin{equation*}
\abs{\mm_{\ell,h}^j(z)}^2
= 1 + \tau^2 \sum_{i=0}^{j-1} \abs{\vv_{\ell,h}^i(z)}^2.
\end{equation*}
Then, both~\eqref{eq:general_constraint} and~\eqref{eq:general_L2_boundedness} 
follow from the equivalence
of the $L^p$-norm of discrete functions with the weighted
$\ell^p$-norm of the vector collecting their nodal values
(with equivalence constants depending only on the shape-regularity of the family of meshes);
see, e.g., \cite[Lemma~3.4]{bartels2015}.
\end{proof}

In the following lemma,
we prove stability of Algorithm~\ref{alg:general}.

\begin{lemma} \label{lem:stability}
Let $\theta_1$, $\theta_2$, $\theta_3$, and $\tau$ satisfy the assumptions of Proposition~\ref{prop:wellposedness_general}
as well as the inequalities
\begin{equation} \label{eq:inequalities_stability}
\theta_1 > 1/2,
\quad
a_{11}a_{22}(2\theta_1-1)^2 > a_{12}^2 (2\theta_2-1)^2,
\quad
\text{and }
\cH^2 \abs{a_0} \abs{2 \theta_3 - 1} \tau < 2.
\end{equation}
Then, there exists a threshold $\tau_0>0$ such that, if $\tau < \tau_0$,
the iterates of Algorithm~\ref{alg:general} satisfy,
for all $j \in \N$,
the stability estimate
\begin{equation} \label{eq:general_stability}
\sum_{\ell=1}^2 \norm[\HH^1(\Omega)]{\mm_{\ell,h}^j}^2
+ \tau \sum_{i=0}^{j-1} \sum_{\ell=1}^2 \norm[\HHH]{\vv_{\ell,h}^i}^2
+ \tau^2 \sum_{i=0}^{j-1} \sum_{\ell=1}^2 \norm{\Grad\vv_{\ell,h}^i}^2
\le C \, \left( 1 + \sum_{\ell=1}^2 \norm[\HH^1(\Omega)]{\mm_{\ell,h}^0}^2 \right).
\end{equation}
The threshold $\tau_0$ depends on $a_0$, $\theta_3$, $\cH$,
and the shape-regularity of the family of meshes,
whereas the constant $C>0$ depends only on $\abs{\Omega}$, $a_{11}$, $a_{12}$, $a_{22}$, $a_0$,
$\theta_1$, $\theta_2$, $\theta_3$,
$\cH$,
and the shape-regularity of the family of meshes.
\end{lemma}

\begin{proof}
Let $j \in \N$.
For all $i=0, \dots, j-1$, we apply Proposition~\ref{prop:energy_law_general}, which yields~\eqref{eq:general_energy_law}.
Summing~\eqref{eq:general_energy_law} over $i=0,\dots,j-1$, we obtain that
\begin{multline*}
\E[\mm_{1,h}^{j},\mm_{2,h}^{j}]
+ \tau \sum_{i=0}^{j-1} \sum_{\ell=1}^2 \norm[\HHH]{\vv_{\ell,h}^i}^2
+ \frac{(2\theta_1-1)}{2} \tau^2 \sum_{i=0}^{j-1} \sum_{\ell=1}^2 a_{\ell\ell} \norm{\Grad\vv_{\ell,h}^i}^2 \\
+ a_{12} (2\theta_2-1) \tau^2 \sum_{i=0}^{j-1} \inner{\Grad\vv_{1,h}^i}{\Grad\vv_{2,h}^i}
- a_0 (2\theta_3-1) \tau^2 \sum_{i=0}^{j-1} \inner{\vv_{1,h}^i}{\vv_{2,h}^i}
=
\E[\mm_{1,h}^0,\mm_{2,h}^0].
\end{multline*}
Using~\eqref{eq:inequalities_stability} and arguing as in the proof of Proposition~\ref{prop:wellposedness_general},
one can show that
\begin{equation*}
\E[\mm_{1,h}^{j},\mm_{2,h}^{j}]
+ \lambda_1 \tau \sum_{i=0}^{j-1} \sum_{\ell=1}^2 \norm[\HHH]{\vv_{\ell,h}^i}^2
+ \lambda_2 \tau^2 \sum_{i=0}^{j-1} \sum_{\ell=1}^2 \norm{\Grad\vv_{\ell,h}^i}^2
\le
\E[\mm_{1,h}^0,\mm_{2,h}^0]
\end{equation*}
for some positive values $\lambda_1 = \lambda_1(a_0, \theta_3)$
and $\lambda_2 = \lambda_2(a_{11}, a_{12}, a_{22}, \theta_1, \theta_2)$.
From~\eqref{eq:assumption_coeff} and Young's inequality,
it follows that
\begin{equation*}
\E[\mm_{1,h}^{j},\mm_{2,h}^{j}]
\ge
\lambda_3 \sum_{\ell=1}^2 \norm{\Grad\mm_{\ell,h}^{j}}^2 
- \frac{\abs{a_0}}{2} \sum_{\ell=1}^2 \norm{\mm_{\ell,h}^{j}}^2
\end{equation*}
for some $\lambda_3 = \lambda_3(a_{11}, a_{12}, a_{22}) > 0$.
Moreover, it holds that
\begin{equation*}
\E[\mm_{1,h}^0,\mm_{2,h}^0]
\le \lambda_4 \sum_{\ell=1}^2 \norm[\HH^1(\Omega)]{\mm_{\ell,h}^0}^2
\end{equation*}
for some $\lambda_4 = \lambda_3(a_{11}, a_{12}, a_{22}, a_0) > 0$.
Altogether, we thus obtain that
\begin{multline} \label{eq:stability_aux1}
\lambda_3 \sum_{\ell=1}^2 \norm{\Grad\mm_{\ell,h}^{j}}^2 
- \frac{\abs{a_0}}{2} \sum_{\ell=1}^2 \norm{\mm_{\ell,h}^{j}}^2
+ \lambda_1 \tau \sum_{i=0}^{j-1} \sum_{\ell=1}^2 \norm[\HHH]{\vv_{\ell,h}^i}^2
+ \lambda_2 \tau^2 \sum_{i=0}^{j-1} \sum_{\ell=1}^2 \norm{\Grad\vv_{\ell,h}^i}^2 \\
\le \lambda_4 \sum_{\ell=1}^2 \norm[\HH^1(\Omega)]{\mm_{\ell,h}^0}^2.
\end{multline}

From Lemma~\ref{lem:projection_free_lemma} and~\eqref{eq:metric},
we deduce that
\begin{equation} \label{eq:stability_aux2}
\abs{a_0} \sum_{\ell=1}^2 \norm{\mm_{\ell,h}^j}^2
\le 2 \cT \abs{a_0} \abs{\Omega}
+ \cT \abs{a_0} \cH \tau^2 \sum_{i=0}^{j-1} \sum_{\ell=1}^2 \norm[\HHH]{\vv_{\ell,h}^i}^2,
\end{equation}
where $\cT>0$ is the constant appearing in~\eqref{eq:general_L2_boundedness}
(which depends only on the shape-regularity of the family of meshes).
Combining~\eqref{eq:stability_aux1} and~\eqref{eq:stability_aux2},
we thus obtain that
\begin{multline*}
\lambda_3 \sum_{\ell=1}^2 \norm{\Grad\mm_{\ell,h}^{j}}^2 
+ \frac{\abs{a_0}}{2} \sum_{\ell=1}^2 \norm{\mm_{\ell,h}^{j}}^2
+ (\lambda_1 - \cT \abs{a_0} \cH \tau) \tau \sum_{i=0}^{j-1} \sum_{\ell=1}^2 \norm[\HHH]{\vv_{\ell,h}^i}^2 \\
+ \lambda_2 \tau^2 \sum_{i=0}^{j-1} \sum_{\ell=1}^2 \norm{\Grad\vv_{\ell,h}^i}^2
\le 2 \cT \abs{a_0} \abs{\Omega}
+ \lambda_4 \sum_{\ell=1}^2 \norm[\HH^1(\Omega)]{\mm_{\ell,h}^0}^2.
\end{multline*}
Hence, if $\tau < \tau_0 := \lambda_1 / (\cT \abs{a_0} \cH)$,
all terms on the left-hand side are nonnegative
and we obtain~\eqref{eq:general_stability},
where the (explicitly computable) constant $C>0$ depends only
on $\abs{\Omega}$, $a_{11}$, $a_{12}$, $a_{22}$, $a_0$,
$\theta_1$, $\theta_2$, $\theta_3$,
$\cH$, and $\cT$.
\end{proof}

In the following proposition,
combining the results we have proved so far,
we establish the main properties of Algorithm~\ref{alg:general}

\begin{proposition} \label{prop:welldefined_general}
Let $\theta_1$, $\theta_2$, $\theta_3$, and $\tau$ satisfy the assumptions of
Lemma~\ref{lem:stability}.
If the time-step size $\tau$ is sufficiently small, then Algorithm~\ref{alg:general} is well defined:
Each iteration is well defined
and the stopping criterion~\eqref{eq:general_stopping} is met in a finite number of iterations.
In particular, the approximate stationary point $(\mm_{1,h} , \mm_{2,h})$ is well defined.
Moreover, for all $\ell=1,2$, it holds that
\begin{equation} \label{eq:constr_error_general}
\norm[L^1(\Omega)]{\interp[\abs{\mm_{\ell,h}}^2]-1}
\le C \tau \left( 1 + \sum_{\ell=1}^2 \norm[\HH^1(\Omega)]{\mm_{\ell,h}^0}^2 \right),
\end{equation}
where
the constant $C>0$ depends only on $\abs{\Omega}$, $a_{11}$, $a_{12}$, $a_{22}$, $a_0$,
$\theta_1$, $\theta_2$, $\theta_3$,
$\cH$,
and the shape-regularity of the family of meshes.
\end{proposition}

\begin{proof}
The well-posedness of each iteration of the algorithm is a consequence of
Proposition~\ref{prop:wellposedness_general}.
Now, let $\tau_0 > 0$ be the threshold guaranteed by Lemma~\ref{lem:stability}.
If $\tau < \tau_0$, then~\eqref{eq:general_stability} holds.
Since the left-hand side of~\eqref{eq:general_stability} is nonnegative
and the right-hand side is independent of $j$,
we deduce that the series 
\begin{equation*}
\sum_{i=0}^{\infty}
\sum_{\ell=1}^2 
\left(\norm[\HHH]{\vv_{\ell,h}^i}^2
+ \tau \norm{\Grad\vv_{\ell,h}^i}^2 \right)
\end{equation*}
is convergent.
It follows that
$\sum_{\ell=1,2}\norm[\HHH]{\vv_{\ell,h}^i}^2
+ \tau \norm{\Grad\vv_{\ell,h}^i}^2 \to 0$ as $i \to \infty$,
which guarantees that
the stopping criterion~\eqref{eq:general_stopping} is satisfied
if $i$ is sufficiently large.
Estimate~\eqref{eq:constr_error_general} is a consequence of~\eqref{eq:general_stability}
and~\eqref{eq:general_constraint} from Lemma~\ref{lem:projection_free_lemma}.
This concludes the proof.
\end{proof}

In the following theorem,
we show the convergence of the sequence generated by Algorithm~\ref{alg:general}.

\begin{theorem}
Let $\theta_1$ and $\theta_2$ satisfy the inequalities
\begin{equation*}
\theta_1>1/2,
\quad
a_{11}a_{22} \theta_1^2 > a_{12}^2 \theta_2^2,
\quad
\text{and}
\quad
a_{11}a_{22} (2\theta_1-1)^2 > a_{12}^2 (2\theta_2-1)^2.
\end{equation*}
Suppose that there exists $c_0>0$,
independent of the discretization parameters $h$, $\tau$, and $\eps$,
such that
\begin{equation} \label{eq:boundedness_initial_general}
\sup_{h>0} \left( \sum_{\ell=1}^2 \norm[\HH^1(\Omega)]{\mm_{\ell,h}^0}^2 \right) \le c_0.
\end{equation}
Suppose that $\tau \to 0$ and $\eps \to 0$ as $h \to 0$.
Then, as $h \to 0$, the sequence of
approximate stationary points
$\{ (\mm_{1,h},\mm_{2,h}) \}_{h>0}$
generated by Algorithm~\ref{alg:general},
upon extraction of a subsequence,
converges weakly in $\HH^1(\Omega) \times \HH^1(\Omega)$
toward a point $(\mm_{1},\mm_{2}) \in \XXX$.
If $a_{12}=0$, the limit $(\mm_{1},\mm_{2})$ is a stationary point
of the energy functional~\eqref{eq:energy}.
\end{theorem}

\begin{proof}
Since $\tau \to 0$, we can assume that it is sufficiently small such that
the algorithm is well defined (cf.\ Proposition~\ref{prop:welldefined_general})
and that the stability estimate~\eqref{eq:general_stability} holds
(cf.\ Lemma~\ref{lem:stability}).
Together with~\eqref{eq:boundedness_initial_general},
it thus follows that the sequence $\{ (\mm_{1,h},\mm_{2,h}) \}_{h>0}$ is uniformly 
bounded in $\HH^1(\Omega) \times \HH^1(\Omega)$.
Hence, there exists $(\mm_1 , \mm_2) \in \HH^1(\Omega) \times \HH^1(\Omega)$
and a (nonrelabeled) weakly convergence subsequence of $\{ (\mm_{1,h},\mm_{2,h}) \}_{h>0}$
such that $(\mm_{1,h},\mm_{2,h}) \weakto (\mm_1 , \mm_2)$ in $\HH^1(\Omega) \times \HH^1(\Omega)$
and $(\mm_{1,h},\mm_{2,h}) \to (\mm_1 , \mm_2)$ in $\LL^2(\Omega) \times \LL^2(\Omega)$.
Combining~\eqref{eq:boundedness_initial_general} with \eqref{eq:constr_error_general},
we see that, for all $\ell = 1,2$, $\norm[L^1(\Omega)]{\interp[\abs{\mm_{\ell,h}}^2]-1} \to 0$ as $h \to 0$.
Hence, applying \cite[Lemma~7.2]{bartels2015}, we obtain that $(\mm_1 , \mm_2) \in \XXX$.

To conclude the proof, it remains to show that,
if $a_{12} = 0$, $(\mm_1 , \mm_2) \in \XXX$ satisfies~\eqref{eq:euler_lagrange1}.
We start with observing that each approximate stationary point $(\mm_{1,h},\mm_{2,h})$
generated by Algorithm~\ref{alg:general}
satisfies the variational formulation
\begin{equation*}
- a_{\ell\ell} \inner{\Grad\mm_{\ell,h}}{\Grad\pphi_{\ell,h}}
+ a_0 \inner{\mm_{3-\ell,h}}{\pphi_{\ell,h}}
= \RRR_{\ell,h}(\pphi_{\ell,h})
\end{equation*}
for all 
$\pphi_{\ell,h} \in \KKK_h[\mm_{\ell,h}]$ and $\ell=1,2$,
where the reminder terms on the right-hand side are given by
\begin{equation*}
\RRR_{\ell,h}(\pphi_{\ell,h}) =
\inner[\HHH]{\vv_{\ell,h}^{i^*}}{\pphi_{\ell,h}}
+ a_{\ell\ell} \theta_1 \tau \inner{\Grad\vv_{\ell,h}^{i^*}}{\Grad\pphi_{\ell,h}}
- a_0  \theta_3 \tau \inner{\vv_{3-\ell,h}^{i^*}}{\pphi_{\ell,h}}
\end{equation*}
and satisfy $\abs{\RRR_h(\pphi_{\ell,h})} \le C \eps \norm[\HH^1(\Omega)]{\pphi_{\ell,h}}$
for all $\pphi_{\ell,h} \in \HHH^1(\Omega)$;
see~\eqref{eq:general1} and \eqref{eq:general_stopping}.
Here, $C>0$ depends only on $a_{11}$, $a_{22}$, $a_0$, and $\abs{\Omega}$.
Note that, since $\eps \to 0$ as $h \to 0$, we have that $\RRR_h \to 0$ in $\HH^1(\Omega)^*$ as $h \to 0$.
Let $\ppsi \in \CC^{\infty}(\overline{\Omega})$.
Choosing the test function $\pphi_{\ell,h} = \interp[\mm_{\ell,h} \times \ppsi] \in \KKK_h[\mm_{\ell,h}]$ in~\eqref{eq:general1}
and passing to the limit as $h \to 0$
(using the available convergence results as in the proof of~\cite[Theorem~7.6]{bartels2015}),
we obtain that
\begin{equation} \label{eq:el_temp}
- a_{\ell\ell} \inner{\Grad\mm_\ell}{\mm_{\ell} \times \Grad \ppsi}
+ a_0 \inner{\mm_{3-\ell}}{\mm_{\ell} \times \ppsi}
= 0
\end{equation}
for all $\ell=1,2$.
Since $\ppsi$ was arbitrary, by density we have that this identity holds for all $\ppsi \in \HH^1(\Omega)$.
Finally, let $\vvphi \in \HH^1(\Omega) \cap \LL^\infty(\Omega)$ be arbitrary.
Choosing $\ppsi = \mm_\ell \times \vvphi$ in~\eqref{eq:el_temp}
and performing simple algebraic manipulations based on the identities
$\aa \times (\bb \times \cc) = (\aa\cdot\cc)\bb - (\aa\cdot\bb)\cc$ (for all $\aa,\bb,\cc \in \R^3$),
$\abs{\mm_\ell}=1$ (a.e.\ in $\Omega$, for all $\ell=1,2$)
and $\partial_i \mm_\ell \cdot \mm_\ell = 0$
(a.e.\ in $\Omega$, for all $i=1,2,3$ and $\ell=1,2$),
we obtain that $(\mm_1 , \mm_2) \in \XXX$ solves~\eqref{eq:euler_lagrange1} for the case $a_{12} = 0$.
This shows that $(\mm_1 , \mm_2)$ is a stationary point of the energy and concludes the proof.
\end{proof}

\subsection{Dynamic problem} \label{sec:proofs_llg}

In this section, we aim to present the proofs of the results concerning Algorithm~\ref{alg:tps}
discussed in Section~\ref{sec:llg}.
However, for the sake of brevity, we omit those of Proposition~\ref{prop:wellposedness}
and Proposition~\ref{prop:energy}, because they can be obtained following line by line
those of Proposition~\ref{prop:wellposedness_general},
Proposition~\ref{prop:energy_law_general},
and Lemma~\ref{lem:stability}.
We focus on the proof of the main convergence result.

\begin{proof}[Proof of Theorem~\ref{thm:main_llg}]
We follow the argument of the seminal paper on the tangent plane scheme~\cite{alouges2008a},
which we adapt in order to take the projection-free update~\eqref{eq:tps2} 
(see also~\cite{ahpprs2014,hpprss2019})
and the different expression of the energy into account.
For the sake of clarity, we split the proof into three steps:
\begin{itemize}
\item \emph{Step~1:} Existence of the limit $(\mm_1, \mm_2) \in L^\infty(0,\infty;\XXX)$.
\end{itemize}
Let $T>0$ be arbitrary.
From the stability estimate~\eqref{eq:llg_stability}
(cf.\ Proposition~\ref{prop:energy}),
which holds uniformly in $h$ and $\tau$ (if $\tau$ is sufficiently small),
it follows that,
for all $\ell=1,2$, the piecewise affine time reconstruction $\mm_{\ell,h\tau}$
and the piecewise constant time reconstructions $\mm_{\ell,h\tau}^\pm$
(defined according to~\eqref{eq:reconstructions})
are both uniformly bounded in
$L^{\infty}(0,\infty;\HH^1(\Omega))$.
Moreover, $\mm_{\ell,h\tau}\vert_{\Omega_T}$
is uniformly bounded in $\HH^1(\Omega_T)$.
By compactness, successive extractions of (nonrelabeled) subsequences
and standard Sobolev embeddings yield
the existence of $\mm_1,\mm_2 \in L^{\infty}(0,\infty;\HH^1(\Omega)) \cap \HH^1(\Omega_T)$
such that, for all $\ell=1,2$, as $h,\tau \to 0$ we have the convergences
$\mm_{\ell,h\tau}\vert_{\Omega_T} \weakto \mm\vert_{\Omega_T}$ in $\HH^1(\Omega_T)$,
$\mm_{\ell,h\tau}\vert_{\Omega_T} \to \mm\vert_{\Omega_T}$ in $\HH^s(\Omega_T)$ for all $s \in (0,1)$,
$\mm_{\ell,h\tau}, \mm_{\ell,h\tau}^{\pm} \weakstarto \mm_\ell$ in $L^{\infty}(0,\infty;\HH^1(\Omega))$,
$\mm_{\ell,h\tau}, \mm_{\ell,h\tau}^{\pm} \weakto \mm_\ell$ in $L^2(0,\infty;\HH^1(\Omega))$,
$\mm_{\ell,h\tau}\vert_{\Omega_T}, \mm_{\ell,h\tau}^{\pm}\vert_{\Omega_T} \to \mm_\ell$ in $L^2(0,T;\HH^s(\Omega))$ for all $s \in (0,1)$,
$\mm_{\ell,h\tau}\vert_{\Omega_T}, \mm_{\ell,h\tau}^{\pm}\vert_{\Omega_T} \to \mm_\ell$ in $\LL^2(\Omega_T)$
and pointwise almost everywhere in $\Omega_T$.
From the projection-free updates~\eqref{eq:tps2},
arguing as in the proof of Lemma~\ref{lem:projection_free_lemma},
we obtain that~\eqref{eq:general_constraint} holds for all $\ell =1,2$ and for all $j \in \N$,
from which it follows (see Step 3 of the proof of~\cite[Proposition~6]{hpprss2019})
that $\abs{\mm_1} = \abs{\mm_2} = 1$ a.e.\ in $\Omega \times (0,\infty)$.
This shows that $(\mm_1, \mm_2) \in L^\infty(0,\infty;\XXX)$.
Finally,
from the stability estimate,
it also follows that, for all $\ell=1,2$,
$\tau \Grad(\mmt_{\ell,h\tau})\vert_{\Omega_T} \to 0$
in $\LL^2(\Omega_T)$
as $h,\tau \to 0$.
\begin{itemize}
\item \emph{Step~2:} If $a_{12}=0$, $(\mm_1, \mm_2)$ satisfies the variational formulation~\eqref{eq:weak:variational}.
\end{itemize}
Let $\vvphi\in\CC^\infty(\overline{\Omega_T})$ be an arbitrary smooth test function.
We consider the smallest integer $j \in \N$ satisfying $T \le j\tau$ and extend $\vvphi$ by zero in $(T,t_j)$.
Let $\ell=1,2$.
For all $i=0,\dots,j-1$, we choose $\pphi_{\ell,h} = \interp[\mm_{\ell,h}^i \times \vvphi(t_i)] \in \KKK_h[\mm_h^i]$
in~\eqref{eq:tps1}, we obtain
\begin{equation*}
\begin{split}
& \alpha_\ell \inner[h]{\vv_{\ell,h}^i}{\interp[\mm_{\ell,h}^i \times \vvphi(t_i)]}
+ \inner[h]{\mm_{\ell,h}^i \times \vv_{\ell,h}^i}{\interp[\mm_{\ell,h}^i \times \vvphi(t_i)]} \\
& \quad + \eta_\ell a_{\ell\ell} \tau \inner{\Grad\vv_{\ell,h}^i}{\Grad\interp[\mm_{\ell,h}^i \times \vvphi(t_i)]} \\
& \qquad =
- \eta_\ell a_{\ell\ell} \inner{\Grad\mm_{\ell,h}^i}{\Grad\interp[\mm_{\ell,h}^i \times \vvphi(t_i)]}
+ \eta_\ell a_0 \inner{\mm_{3-\ell,h}^i}{\interp[\mm_{\ell,h}^i \times \vvphi(t_i)]},
\end{split}
\end{equation*}
Due to the properties of the mass-lumped scalar product,
we can remove the nodal interpolant from the first two terms on the left-hand side
without altering the value of the integrals.
Multiplication by $\tau$ and summation over $i=0,\dots,j-1$ then yield
\begin{equation*}
\begin{split}
& \alpha_\ell \int_0^{t_j}\inner[h]{\mmt_{\ell,h\tau}(t)}{\mm_{\ell,h\tau}^-(t) \times \vvphi_\tau^-(t)} \dt \\
& \quad
+ \int_0^{t_j} \inner[h]{\mm_{\ell,h\tau}^-(t) \times \mmt_{\ell,h\tau}(t)}{\mm_{\ell,h\tau}^-(t) \times \vvphi_\tau^-(t)]} \dt \\
& \qquad
+ \eta_\ell a_{\ell\ell} \tau \int_0^{t_j} \inner{\Grad\mmt_{\ell,h\tau}(t)}{\Grad\interp[\mm_{\ell,h\tau}^-(t) \times \vvphi_\tau^-(t)]} \dt \\
& \qquad\quad =
- \eta_\ell a_{\ell\ell} \int_0^{t_j} \inner{\Grad\vvphi_\tau^-(t)}{\Grad\interp[\mm_{\ell,h\tau}^-(t) \times \vvphi_\tau^-(t)]} \dt \\
& \qquad \qquad + \eta_\ell a_0 \int_0^{t_j} \inner{\mm_{3-\ell,h\tau}^-(t)}{\interp[\mm_{\ell,h\tau}^-(t) \times \vvphi_\tau^-(t)]} \dt,
\end{split}
\end{equation*}
where we note that we have rewritten the equation in terms of the time reconstructions~\eqref{eq:reconstructions}.
Using~\eqref{eq:h-inner-product} and the approximation properties of the nodal interpolant,
in all integrals we substitute the mass-lumped inner products by $L^2$-products
and remove the nodal interpolant (see~\cite{alouges2008a}).
Moreover, exploiting the fact that the integrands are all uniformly bounded,
we modify the domain in integration in time from $(0,t_j)$ to $(0,T)$.
All these actions generate an error which goes to zero in the limit as $h,\tau \to 0$.
In particular, we obtain
\begin{equation*}
\begin{split}
& \alpha_\ell \int_0^T \inner{\mmt_{\ell,h\tau}(t)}{\mm_{\ell,h\tau}^-(t) \times \vvphi_\tau^-(t)} \dt \\
& \quad + \int_0^T \inner{\mm_{\ell,h\tau}^-(t) \times \mmt_{\ell,h\tau}(t)}{\mm_{\ell,h\tau}^-(t) \times \vvphi_\tau^-(t)]} \dt \\
& \qquad
+ \eta_\ell a_{\ell\ell} \tau \int_0^T \inner{\Grad\mmt_{\ell,h\tau}(t)}{\Grad[(\mm_{\ell,h\tau}^-(t) \times \vvphi_\tau^-(t)]} \dt \\
& \quad\qquad =
- \eta_\ell a_{\ell\ell} \int_0^T \inner{\Grad\vvphi_\tau^-(t)}{\Grad[\mm_{\ell,h\tau}^-(t) \times \vvphi_\tau^-(t)]} \dt \\
& \qquad\qquad
+ \eta_\ell a_0 \int_0^T \inner{\mm_{3-\ell,h\tau}^-(t)}{\mm_{\ell,h\tau}^-(t) \times \vvphi_\tau^-(t)} \dt
+ o(1).
\end{split}
\end{equation*}
Using the convergence results available from Step~1,
we can pass this formulation to the limit as $h,\tau \to 0$
and obtain that the last term on the left-hand side goes to zero,
whereas all other terms converge toward the corresponding ones in~\eqref{eq:weak:variational}.
For the details of the argument, we refer to~\cite{alouges2008a} for all terms
but the second one on the left-hand side, which, due to the omission of the nodal projection from~\eqref{eq:tps2},
requires a more careful treatment (see Step~2 of the proof of~\cite[Theorem~1]{hpprss2019}).
This shows that, for all $\ell=1,2$, $\mm_\ell$ satisfies~\eqref{eq:weak:variational}
for all $\vvphi \in \CC^\infty(\overline{\Omega_T})$.
By density, the result then holds for all $\vvphi \in \HH^1(\Omega_T)$.
\begin{itemize}
\item \emph{Step~3:} $(\mm_1, \mm_2)$ satisfies the energy inequality~\eqref{eq:weak:energy}.
\end{itemize}
We start from the discrete energy law~\eqref{eq:llg_energy_law} established in Proposition~\ref{prop:energy}.
Using~\eqref{eq:assumption_coeff} and a combination of Cauchy--Schwarz' an Young's inequalities, we obtain that
\begin{multline*}
\E[\mm_{1,h}^j,\mm_{2,h}^j]
+ \tau \sum_{i=0}^{j-1} \sum_{\ell=1}^2 \left( \frac{\alpha_\ell}{\eta_\ell} - \frac{\abs{a_0}\tau}{2}\right) \norm[h]{\vv_{\ell,h}^i}^2
+ \lambda \tau^2 \sum_{i=0}^{j-1} \sum_{\ell=1}^2 \norm{\Grad\vv_{\ell,h}^{i}}^2
\le
\E[\mm_{1,h}^0,\mm_{2,h}^0],
\end{multline*}
where $\lambda>0$ is the minimum eigenvalue
of the 2-by-2 matrix $\begin{pmatrix} a_{11} & a_{12} \\ a_{12} & a_{22} \end{pmatrix}$.
The last term on the left-hand side is nonnegative and can be omitted.
Rewriting the inequality in terms of the time reconstructions~\eqref{eq:reconstructions},
we get
\begin{equation*}
\E[\mm_{1,h\tau}^+(T),\mm_{2,h\tau}^+(T)]
+ \sum_{\ell=1}^2 \left( \frac{\alpha_\ell}{\eta_\ell} - \frac{\abs{a_0}\tau}{2}\right) \int_0^T \norm[h]{\mmt_{\ell,h\tau}(t)}^2 \, \dt
\le
\E[\mm_{1,h\tau}^-(0),\mm_{2,h\tau}^-(0)].
\end{equation*}
Passing to the limit as $h,\tau \to 0$,
using the convergence results available from Step~1,
standard lower semicontinuity arguments yield~\eqref{eq:weak:energy}.
This concludes the proof.
\end{proof}

\section{Acknowledgment}

MR is a member of the `Gruppo Nazionale per il Calcolo Scientifico (GNCS)'
of the Italian `Istituto Nazionale di Alta Matematica (INdAM)'.
Part of the work on this paper was undertaken when the authors were visiting
the 
Hausdorff Research Institute for Mathematics
of the University of Bonn
during the Trimester Program on
\emph{Mathematics for Complex Materials},
funded by the German Research Foundation (DFG)
under Germany's Excellence Strategy -- EXC-2047/1-- 390685813.
The kind hospitality of the institute is thankfully acknowledged.

\bibliographystyle{habbrv}
\bibliography{ref}

\appendix

\section{The equations in physical units} \label{sec:nondimensional}

In this appendix,
for the convenience of all interdisciplinary readers,
we present the model in physical units
(used for physical investigations, e.g., in~\cite{ne2015,mlp2016,pkcatsf2019,mra2019,sztfm2020,spkcf2020,tslggcaf2020,cstcf2021})
and show how to obtain from it the dimensionless setting
described in Section~\ref{sec:model} and analyzed in the paper.
By doing this, we also justify the setup and the choice of the material parameters
in the numerical experiments presented in the work.

\subsection{Nondimensionalization} \label{sec:physical}

Let $\Omega \subset \R^3$ be the volume occupied by an AFM or FiM material.
Let the vector field $\MM : \Omega \to \R^3$ denote the total magnetization of the sample (in \si{\ampere\per\meter}).
The total magnetization can be decomposed as
$\MM = \MM_1 + \MM_2$,
where $\MM_1 , \MM_2 : \Omega \to \R^3$,
the magnetization vectors of two sublattices (in \si{\ampere\per\meter}),
satisfy the constraints $\abs{\MM_1} = \Msl{1}$ and $\abs{\MM_2} = \Msl{2}$.
The constants $\Msl{1},\Msl{2} > 0$ are the sublattice saturation magnetizations  (in \si{\ampere\per\meter}).
Let $\mm_1,\mm_2 : \Omega \to \sphere$ be
the dimensionless unit-length vector fields
$\mm_1 = \MM_1/\Msl{1}$ and $\mm_2 = \MM_2/\Msl{2}$.
The total Gibbs free energy (in \si{\joule}) of the system
(assumed, for simplicity, to include only exchange contributions in this section)
reads as
\begin{equation} \label{eq:energy_appendix}
\begin{split}
\E [\mm_1,\mm_2]
& = \E_{\mathrm{ex}} [\mm_1,\mm_2] \\
& = \sum_{\ell=1}^2 A_{\ell\ell} \int_{\Omega} \abs{\Grad\mm_\ell}^2
+ A_{12} \int_{\Omega} \Grad\mm_1 : \Grad\mm_2
- \frac{4 A_0}{a^2} \int_{\Omega} \mm_1 \cdot \mm_2,
\end{split}
\end{equation}
where the exchange constants $A_{11}, A_{22}>0$ and $A_{12} ,A_0 \in \R$ are in \si{\joule\per\meter},
whereas $a>0$ is the lattice constant (in \si{\meter}).
The first contribution in~\eqref{eq:energy_appendix} is called
\emph{inhomogeneous intralattice exchange}
and models the classical ferromagnetic exchange for $\mm_1$ and $\mm_2$.
The second term is called
\emph{inhomogeneous interlattice exchange},
which arises from a nearest-neighbor approximation of the exchange interaction
between spins.
The last contribution is called
\emph{homogeneous interlattice exchange}
and takes the local interaction between $\mm_1$ and $\mm_2$ into account.

The dynamics of $\mm_1$ and $\mm_2$ is governed by a coupled system of two LLG equations
\begin{equation} \label{eq:LLG}
\mmt_\ell = - \gamma_\ell \, \mm_\ell \times \Heff{\ell}[\mm_1,\mm_2] + \alpha_\ell \, \mm_\ell \times \mmt_\ell
\quad
\text{for } \ell = 1,2,
\end{equation}
where $\gamma_\ell > 0$ (in \si{\meter\per\ampere\per\second}) and $\alpha_\ell > 0$ (dimensionless)
are the sublattice rescaled gyromagnetic ratios and Gilbert damping parameters, respectively.
In~\eqref{eq:LLG}, the effective fields $\Heff{\ell}[\mm_1,\mm_2]$ (in \si{\ampere\per\meter})
are equal, up to a negative multiplicative constant, to the functional derivatives of the total energy with respect to $\mm_\ell$, i.e.,
\begin{equation*}
\Heff{\ell}[\mm_1,\mm_2] = - \frac{1}{\mu_0 \Msl{\ell}} \frac{\E [\mm_1,\mm_2]}{\delta m_\ell},
\end{equation*}
where
$\mu_0$ is the vacuum permeability (in \si{\newton\per\ampere\squared}).
Assuming no flux boundary conditions, the strong form of the resulting effective fields reads as
\begin{equation*}
\Heff{\ell}[\mm_1 , \mm_2]
= \frac{2A_{\ell\ell}}{\mu_0 \Msl{\ell}} \, \Lapl \mm_\ell
+ \frac{A_{12}}{\mu_0 \Msl{\ell}} \, \Lapl \mm_{3 - \ell}
+ \frac{4 A_0}{\mu_0 \Msl{\ell} \, a^2} \, \mm_{3 - \ell}.
\end{equation*}
We now start the nondimensionalization.
Let $\Ms>0$ and $\gamma_0>0$ be
some reference saturation magnetization (in \si{\ampere\per\meter})
and rescaled gyromagnetic ratio (in \si{\meter\per\ampere\per\second}),
respectively.
For all $\ell = 1,2$, define the positive dimensionless parameters $\eta_{\mathrm{s},\ell} := \Msl{\ell} / \Ms$
and $\eta_{\mathrm{\gamma},\ell} := \gamma_\ell / \gamma_0$.
The dimensionless total magnetization is given by
$\mm = \MM/\Ms = \eta_{\mathrm{s},1} \mm_1 + \eta_{\mathrm{s},2} \mm_2$.

Let $L>0$ is some intrinsic length of the problem.
We rescale the space and time variables to obtain the dimensionless variables
$x' = x/L$ and $t' = \gamma_0 \Ms \, t$.
Accordingly, we rescale also the domain $\Omega' = \Omega/L$.
We consider the rescaled unit-length vector fields $\mm_\ell'(x',t') = \mm_\ell(L x' , t' / (\gamma_0\Ms))$
($\ell=1,2$)
and the rescaled total magnetization
$\mm'(x',t') = \mm(L x' , t' / (\gamma_0\Ms))$.
Moreover, we rescale the energy as $\E'[\mm_1',\mm_2'] = \E[\mm_1,\mm_2] / (\mu_0 \Ms^2 L^3)$,
which yields the expression
\begin{multline*}
\E' [\mm_1',\mm_2']
=
\E_{\mathrm{ex}}' [\mm_1',\mm_2'] \\
= \frac{1}{2}\sum_{\ell=1}^2 \frac{2A_{\ell\ell}}{\mu_0 \Ms^2 L^2} \int_{\Omega'} \abs{\Grad'\mm_\ell'}^2
+ \frac{A_{12}}{\mu_0 \Ms^2 L^2} \int_{\Omega'} \Grad'\mm_1' : \Grad'\mm_2'
- \frac{4 A_0}{\mu_0 \Ms^2 a^2} \int_{\Omega'} \mm_1' \cdot \mm_2'.
\end{multline*}
Defining the dimensionless coefficients
$a_{\ell\ell} = 2A_{\ell\ell}/(\mu_0 \Ms^2 L^2)>0$ ($\ell=1,2$),
$a_{12} = A_{12}/(\mu_0 \Ms^2 L^2) \in \R$,
and $a_0 = 4 A_0/(\mu_0 \Ms^2 a^2) \in \R$,
and omitting all `primes' for simplicity,
we obtain the dimensionless energy functional~\eqref{eq:energy} of Section~\ref{sec:model}.
By construction, the dimensionless rescaled effective fields defined in \eqref{eq:heff}
are related to the ones in physical units according to the relation
\begin{equation*}
\heff{\ell}[\mm_1',\mm_2']
\stackrel{\eqref{eq:heff}}{=} - \frac{\delta\E'[\mm_1',\mm_2']}{\delta \mm_{\ell}'}
= \frac{\eta_{\mathrm{s},\ell}^2}{\Msl{\ell}} \, \Heff{\ell}[\mm_1,\mm_2]
\quad
\text{for all } \ell = 1,2.
\end{equation*}
Rescaling the LLG equations in~\eqref{eq:LLG} according to the above change of variables
and introducing all dimensionless quantities,
we obtain
\begin{equation*}
\partial_{t'} \mm_\ell'
= - \frac{\eta_{\gamma,\ell}}{\eta_{\mathrm{s},\ell}} \, \mm_\ell' \times \heff{\ell}[\mm_1',\mm_2'] + \alpha_\ell \, \mm_\ell' \times \partial_{t'}\mm_\ell'
\quad
\text{for all } \ell = 1,2,
\end{equation*}
Defining the dimensionless parameter $\eta_\ell := \eta_{\gamma,\ell} / \eta_{\mathrm{s},\ell} > 0$
and omitting all `primes',
we obtain the dimensionless system~\eqref{eq:llg} of LLG equations of Section~\ref{sec:model}.

\subsection{Lower-order energy contributions} \label{sec:lower}

In practically relevant simulations,
to be able to describe complex physical processes involving AFM and FiM materials,
more energy contributions (in addition to the exchange ones) need to be taken into account in~\eqref{eq:energy_appendix}:
\begin{itemize}
\item
The \emph{magnetocrystalline anisotropy energy}
incorporates the existence of preferred directions of alignment for the fields.
In the uniaxial case, it reads as
\begin{equation*}
\E_{\mathrm{ani}} [\mm_1,\mm_2]
=
K_1 \int_{\Omega} [1 - (\aa_1\cdot\mm_1)^2]
+ K_2 \int_{\Omega} [1 - (\aa_2\cdot\mm_2)^2],
\end{equation*}
where $K_1,K_2 > 0$ are physical constants (in \si{\joule\per\meter\cubed}),
whereas
$\aa_1,\aa_2 \in\sphere$ are the so-called easy axes of the material
(usually it holds that $K_1=K_2$ and $\aa_1=\aa_2$).

\item
The \emph{Dzyaloshinskii--Moriya interaction} is used to incorporate chiral effects into the model.
Its general expression for AFM and FiM materials is given by
\begin{equation*}
\E_{\mathrm{DMI}} [\mm_1,\mm_2]
=
\int_{\Omega} \DD_1 : (\Grad\mm_1 \times \mm_1)
+ \int_{\Omega} \DD_2 : (\Grad\mm_2 \times \mm_2),
\end{equation*}
where $\DD_1,\DD_2 \in \R^{3 \times 3}$ are the so-called spiralization tensors
(with coefficients in \si{\joule\per\meter\squared}),
whereas, for $\ell=1,2$,
$\Grad\mm_\ell \times \mm_\ell$ denotes the matrix with columns
$\partial_j \mm \times \mm$ for $j=1,2,3$
(again, usually it holds that $\DD_1=\DD_2$).

\item
The \emph{Zeeman energy} models the interaction of the total magnetization with
an applied external field (assumed to be magnetization-independent)
and reads as
\begin{equation*}
\E_{\mathrm{ext}} [\mm_1,\mm_2]
=
- \mu_0 \int_{\Omega} \Hext\cdot(\Msl{1}\mm_1 + \Msl{2}\mm_2),
\end{equation*}
where $\Hext \in \R^3$ denotes an applied external field (in \si{\ampere\per\meter}).

\item
The \emph{magnetostatic energy}
can be understood as the energy associated with the interaction
of the total magnetization with the stray field $\Hstray \in \R^3$,
which solves the magnetostatic Maxwell equations
\begin{equation*}
\div\Hstray = - \div[\chi_\Omega (\Msl{1}\mm_1 + \Msl{2}\mm_2)]
\quad
\text{and}
\quad
\curl\Hstray = \0
\quad
\text{in } \R^3.
\end{equation*}
The energy contribution is given by
\begin{equation*}
\E_{\mathrm{ext}} [\mm_1,\mm_2]
=
- \frac{\mu_0}{2} \int_{\Omega} \Hext\cdot(\Msl{1}\mm_1 + \Msl{2}\mm_2),
\end{equation*}
where $\chi_\Omega : \R^3 \to \{0,1\}$ denotes the indicator function of the domain $\Omega$.
\end{itemize}
Note that in all the above energy contributions the two fields are decoupled
(for the magnetostatic energy, this is a consequence of the fact that the operator
mapping the total magnetization to the solution of the
magnetostatic Maxwell equations is linear).
Hence, even in the presence of the above contributions,
the system of Euler--Lagrange equations associated with the minimization problem
and the system of LLG equations are only exchange-coupled.

In the numerical experiments of the work (see Sections~\ref{sec:numerics_stationary} and~\ref{sec:numerics_llg}),
we considered dimensionless forms of magnetocrystalline anisotropy energy,
Dzyaloshinskii--Moriya interaction and Zeeman energy, namely
\begin{align*}
\E_{\mathrm{ani}} [\mm_1,\mm_2]
&=
\frac{q_1^2}{2} \int_{\Omega} [1 - (\aa_1\cdot\mm_1)^2]
+ \frac{q_2^2}{2} \int_{\Omega} [1 - (\aa_2\cdot\mm_2)^2],\\
\E_{\mathrm{DMI}} [\mm_1,\mm_2]
&=
\int_{\Omega} \widehat\DD_1 : (\Grad\mm_1 \times \mm_1)
+ \int_{\Omega} \widehat\DD_2 : (\Grad\mm_2 \times \mm_2),\\
\E_{\mathrm{ext}} [\mm_1,\mm_2]
&=
- \int_{\Omega} \hext\cdot(\eta_{\mathrm{s},1} \mm_1 + \eta_{\mathrm{s},2} \mm_2)
=
- \int_{\Omega} \hext\cdot\mm.
\end{align*}
In these expressions, which can be obtained rescaling the energy contributions as described in the previous section,
the dimensionless parameters are related to the physical ones via the relationships
$q_\ell = \sqrt{2 K_\ell/(\mu_0 \Ms^2)}$,
$\widehat \DD_\ell = \DD_\ell/(\mu_0 \Ms^2 L)$ ($\ell=1,2$),
and $\hext = \Hext/\Ms$.

To conclude, we note that for AFM and FiM materials, differently from what happens for FM materials,
the Zeeman and the magnetostatic energies are usually of limited physical importance,
because they depend on the total magnetization of the sample, which is in general very small.

\end{document}